\documentclass[11pt,a4paper,reqno,twoside]{amsart}

\usepackage{amsmath,amsthm,amsfonts,amssymb,euscript,cite, mathrsfs,mathtools}
\usepackage{xcolor} %
\usepackage{hyperref}
\hypersetup{colorlinks,linkcolor={blue},citecolor={blue},urlcolor={red}}  
\usepackage{enumitem} %
\usepackage{ulem}
\usepackage[capitalise]{cleveref}
\usepackage{verbatim}
\crefname{equation}{}{}

\usepackage{fullpage} %

\usepackage{seqsplit}
\usepackage{fancyhdr}

\definecolor{green}{rgb}{0,0.8,0} %

\definecolor{babypink}{rgb}{0.96,0.76,0.76}

\newcommand{\ang}{{\not\negmedspace\partial }}

\newcommand{\rt}{{\mathbb R^3}}

\newcommand{\pa}{\partial}

\newcommand{\gab}{g^{\alpha\beta}}
\newcommand{\hab}{h^{\alpha\beta}}

\newcommand{\pab}{\partial_\beta}
\newcommand{\paa}{\partial_\alpha}

\newcommand{\pat}{\partial_t}

\newcommand{\ti}{\tilde}

\newcommand{\la}{\langle}
\newcommand{\ra}{\rangle}

\newcommand{\ls}{\lesssim}

\newcommand{\inv}{^{-1}}

\newcommand{\f}{\frac}

\newcommand{\iy}{\infty}
\renewcommand{\S}{{\mathbb S}}

\newcommand{\crt}{{C^{ R}_{ T}}}
\newcommand{\cut}{{C^{ U}_{ T}}}

\newcommand{\inte}{{C^{\rm{int}}_{ T}}}
\newcommand{\tinte}{{\ti C^{\rm{int}}_T}}

\newcommand{\lolt}{{L^1L^{2}}}

\newcommand{\lt}{{L^2}}
\newcommand{\p}{\phi}

\renewcommand{\pm}{\phi_{\le m}}

\newcommand{\pmn}{\p_{\le m+n}}
\newcommand{\nm}{\la t-r\ra} %
\newcommand{\jr}{\la r\ra} %
\newcommand{\ju}{\la u\ra}

\newcommand{\jt}{\la t\ra}
\newcommand{\jv}{\la v\ra}
\newcommand{\js}{\la s\ra}
\newcommand{\jrho}{\la\rho\ra}

\newcommand{\lr}[1]{\left( #1 \right)}

\newcommand{\norm}[1]{\left\lVert#1\right\rVert}

\let\arXiv\arxiv
\def\doi#1{ {\href{http://dx.doi.org/#1}
   {{\mdseries\ttfamily DOI}}}}

\newtheorem{theorem}{Theorem}[section]

\newtheorem{lemma}[theorem]{Lemma}
\newtheorem{proposition}[theorem]{Proposition}
\theoremstyle{definition}

\newtheorem{definition}[theorem]{Definition}

\theoremstyle{remark}
\newtheorem{remark}[theorem]{Remark}
\numberwithin{equation}{section}

\newcommand{\x}{\alpha}
\newcommand{\xb}{\beta}

\newcommand{\eps}{\epsilon}

\newcommand{\xo}{\omega}
	\newcommand{\xO}{\Omega}

\newcommand{\xs}{\sigma}

\newcommand{\N}{{\mathbb N}}

\newcommand{\R}{\mathbb R}

\newcommand{\cal}{\mathcal }

\newcommand{\calD}{\mathcal D}
\newcommand{\calE}{\mathcal E}

\newcommand{\calN}{\mathcal N}

\newcommand{\calR}{\mathcal R}

\newcommand{\bo}{\Box}

\newcommand{\tpa}{\bar\pa}

\usepackage{microtype} %
\begin{document}

\title[Global existence and pointwise decay for null condition]{Global existence and pointwise decay for nonlinear waves under the null condition} 

\author[S. Looi]{Shi Zhuo Looi}
\address{The Division of Physics, Mathematics, and Astronomy, California Institute of Technology, Pasadena, CA 91125}
\email{looi@caltech.edu}
\author[M. Tohaneanu]{Mihai Tohaneanu}
\address{Department of Mathematics, University of Kentucky, Lexington, KY 40506}
\email{mihai.tohaneanu@uky.edu}

\begin{abstract}
This paper proves global existence and sharp pointwise decay for solutions to nonlinear wave equations satisfying the semilinear null condition, on a class of three-dimensional, asymptotically flat, and notably, non-stationary spacetimes. We consider nonlinearities satisfying a generalized null condition which does not necessarily retain its structure when commuted with vector fields. For sufficiently small initial data, and under the assumption that the underlying linear operator satisfies an integrated local energy decay estimate, we prove that solutions exist for all time and we establish sharp pointwise decay estimates for the solution $\phi$ and its vector-fields. The solution itself decays as $|\phi(t,x)| \lesssim \langle t+r \rangle^{-1} \langle t-r \rangle^{-1}$. This rate matches that of the nonlinear equation on a flat background. This rate is sharp, as this behavior holds already for certain time-dependent perturbations of the classical null form on Minkowski space, which we specify. 
\end{abstract}

\maketitle

\section{Introduction}
In this paper, we study nonlinear wave equations subject to the generalized null condition in various space-times. Our objectives are twofold: first, to establish global existence of solutions, and second, to derive sharp pointwise decay (upper bounds on the absolute value of the solution) estimates in a robust manner. The nonlinearities addressed here extend beyond the classical null condition, and are a generalization of the classical null condition. Unlike the classical case, the null structure is not preserved when commuting with derivatives or vector fields, which adds an extra challenge that we address.

We also emphasize that we consider dynamical (time-dependent) perturbations of $\Box$, where $\Box := - \partial_t^2 + \Delta$ is the standard wave operator on $\R\times \R^3$. Not many previous works have considered time-dependent perturbations of $\Box$. We expect that the rate of pointwise decay we obtain for the solution is sharp because the rate is expected to be sharp in the flat case, i.e., the special case of $P = \Box$. Our method not only provides pointwise decay rates for the solution but also, in most cases, gives improved decay rates for its derivatives and vector fields.

The paper is structured as follows. Section 1 introduces our result, and some of the history of the problem. Section 2 contains some notation, a discussion of local energy estimates, and the rigorous statement of our main theorem. Sections 3 and 4 contain the proof of global existence. Sections 5, 6 and 7 are dedicated to the sharp pointwise bounds.  

\subsection{Statement of the result}

The goal of this paper is to study wave equations of the form $P\phi= F(\partial\p)$, where the operator $P$ is defined by \begin{equation}\label{P def}
P := \paa\gab(t,x)\pab + g^\xo(t,x) \Delta_\xo + B^\x(t,x)\paa + V(t,x).
\end{equation} 
$\Delta_\xo$ denotes the Laplace operator on the unit sphere, $\x,\xb$ range across $0, \dots, 3$, and $g^{\alpha\beta}$, $g^\xo$, $B^\x$ and $V$ denote smooth functions.  We assume that $P$ is hyperbolic, asymptotically flat, and that the linear evolution satisfies strong local energy decay (and thus the Hamiltonian flow must be nontrapping); the precise conditions on the potential $V$, the coefficients $B,g^\xo$ and the Lorentzian metric $g$ are given in the main result, \cref*{thm:main}. We permit time-dependent coefficients, and also allow for large perturbations of $\Box$ within compact spatial regions. 

The non-linearity $F$ is a quadratic expression in $\pa\p$ that generalizes the classical null condition of Klainerman \cite{K2}. For a semilinear wave equation with constant coefficients, $\Box \p = F(\pa\p)$, the classical null condition amounts to requiring that 
\begin{equation}\label{null}
F(\pa\p)=\lambda ((\pa_t\p)^2 - |\pa_x\p|^2)),
\end{equation} 
where $\lambda$ is any constant. One thus considers
\begin{equation} \label{eq:classicalproblem}
\Box\p= \lambda((\pa_t\p)^2 - |\pa_x\p|^2)), \ \ (\phi(0,x), \pat\phi(0,x)) =(\p_0(x),\p_1(x)), 
\end{equation} 
 For the equation \cref*{eq:classicalproblem}, given a pair of smooth and compactly supported functions $f$ and $g$, there exists a small number $\varepsilon_0 > 0$ depending on this pair of functions, such that if the initial data satisfy
\[
(\phi(0,x), \pat\phi(0,x)) = (\varepsilon f, \varepsilon g),
\]
then the solution is global in time whenever $\varepsilon \leq \varepsilon_0$. 

We shall prove global existence and pointwise decay for \cref*{eq:problem}, which is a generalized version of \cref*{eq:classicalproblem}. The nonlinearity in \cref*{eq:problem} (see \cref*{Qass}) is a generalized version of the nonlinearity \cref*{null} defined in the constant-coefficient case, and in \cref*{eq:problem}, the wave operator $P$ is a significant generalization of $\Box$, see \cref*{thm:main}.  We invite the reader to consider our nonlinear condition as a set of sufficient criteria for variable-coefficient quadratic nonlinearities to not deviate too far from, in terms of nonlinear estimates, the behavior observed under the classical constant-coefficient quadratic null condition for systems of equations. However, the nonlinearities we will consider may depart significantly from the classical case. One such example is any $F$ that outside of the unit ball takes the form
\[(\partial_t\phi)^2 + \frac{x^i}{r} \partial_{x^i}\phi \, \partial_t\phi = \partial_t\phi(\partial_t\phi + \partial_r\phi).\] 
 Another example would be a nonlinearity that, outside the unit ball, takes the form $\partial_\omega\phi \partial \phi$, where $\partial_\omega$ denotes any of the three possible angular derivatives defined in \cref*{def:angular}, and $\partial \phi$ denotes any derivative of $\phi$.
We also allow for general coefficients that are functions in the class $S^Z(1)$ (see \cref*{SZdef}). Thus, for example, $$\frac{t^2+r^2}{t^2 + r^2+1} \left((\partial_t\phi)^2 + \frac{x^i}{r} \partial_{x^i}\phi \, \partial_t\phi \right)\chi_{|x| > 1}(t,x)$$ would be allowed, where $\chi_{|x| > 1}$ is a smooth cutoff function supported away from $x = 0$.

Our main theorem asserts that if the solution to the linear wave equation $P\phi = F$ satisfies a strong local energy decay estimate, then the problem described by \cref*{eq:problem} with sufficiently small initial data admits a unique global solution. Moreover, this solution and the vector fields $Z$ applied to it have global space-time pointwise decay rates of
\[
\langle t - r \rangle^{-1} \langle t + r \rangle^{-1}.
\] The decay rate we obtain matches the one established by Christodoulou \cite{Ch} for the case $P = \Box$ using the conformal method. We conjecture that this rate is sharp for the level of generality considered with the wave operator $P$. For the precise formulation of our main theorem, see \cref*{thm:main}.

The proof strategy begins with writing the solution as $\phi = w + \psi$ where $w$ solves $\Box w=0$ with the same initial data as $\phi$. Given appropriate initial data, the function $w$ already has the decay rate stated in Theorem~\ref{thm:main}. Then we reframe the equation for $\psi$ as 
$$\Box\psi = (\Box-P)\phi + F(\partial\phi),$$
where $\Box$ is the d'Alembertian operator, and $\psi(0,x) = 0, \partial_t \psi(0,x) = 0$. This formulation allows us to leverage results for the inhomogeneous wave equation
$$\Box u = f$$
which we prove in Section~\ref{sec:preliminaries}. We apply these results with $u = \psi$ and $f = (\Box-P)\phi + F(\pa \phi)$.
We utilize $L^\infty$-$L^\infty$ estimates for the fundamental solution of $\Box$ on $\R \times \R^3$ in order to carry out an iteration process which gradually improves the decay rate of the solution $\psi$. These $L^\infty$-$L^\infty$ estimates show that for an inhomogeneity $f$ with a given space-time decay rate, the solution $u$ exhibits improved decay properties. Standard tools like those found in Section~\ref{sec:pointwise} give a decay rate for the solution, and our proof strategy is to iteratively improve this decay rate. The results we prove which show the improved decay rates of derivatives of the solution in Section~\ref{sec:preliminaries} are some of the more challenging parts of this iteration process.

\subsection{History}

The semilinear wave equation in $\R^{1+3}$
\begin{equation}\label{Mnull}
\Box \phi = Q[\pa\p,\pa\p], \qquad \phi |_{t=0} = \phi_0, \qquad \pa_t \phi |_{t=0} = \phi_1
\end{equation}
for small initial data has been studied extensively. 

A first natural (and nontrivial) question to ask is whether there is a unique global solution. It is known that the answer is no, at least for general nonlinearities; indeed, solutions blows up in finite time if $Q[\pa\p,\pa\p] = (\pa_t\p)^2$, see \cite{J}. On the other hand, if the nonlinearity satisfies the null condition \cref*{null}, first identified by Klainerman \cite{K1}, it was shown independently in \cite{Ch} and \cite{K2} that the solution exists globally. This result was extended to quasilinear systems with multiple speeds, as well as the case of exterior domains; see, for instance, \cite{MNS}, \cite{MS1}, \cite{MS2}, \cite{H}, \cite{S}, \cite{KS}, \cite{Al}, \cite{Lin1}, \cite{ST}, \cite{FM}, as well as to systems satisfying the weak null condition, including Einstein's Equations, see \cite{LR1}, \cite{LR2}, \cite{LR3}.

There have also been many works for small data in the variable coefficient case. Almost global existence for nontrapping metrics was shown in \cite{BH}, \cite{SW}. Global existence for stationary, small perturbations of Minkowski was shown in \cite{WY}, and for nonstationary, compactly supported perturbations in \cite{Y}.  In the context of black holes, global existence was shown in \cite{Luk} for Kerr space-times with small angular momentum, and in \cite{AAG} for the Reissner-Nordstr\"om backgrounds. 

A second key question, under the assumption of global existence, is to analyze the asymptotic behavior of solutions as $t \to \infty$. A natural approach to this is to estimate the rate of pointwise decay of the solution. Moreover, as part of this investigation, weaker decay rates are not only achievable—and indeed attained—but also indispensable for the standard approach to establishing global existence. However, in order to completely understand the behavior of solutions, one would like to obtain sharp rates of decay. For linear wave equations on Schwarzschild and Kerr backgrounds, sharp pointwise rates of $t^{-3}$ were conjectured by Price \cite{Pri} and recently proved rigorously by several authors \cite{DSS}, \cite{Tat}, \cite{MTT}. Recently, these sharp decay rates were also used to obtain precise  asymptotics for linear wave equations on subextremal Kerr black holes \cite{AAG2}, \cite{Hintz}. 

For nonlinear problems, sharp rates of decay have been established for power nonlinearities \cite{Toh}, \cite{Looi2}, \cite{Looi}. See also the very recent \cite{LukOh,LX} for sharp pointwise bounds and asymptotics, given certain assumptions, for a variety of nonlinearities. The work \cite{LukOh} shows, among a large number of other results, that the decay rate obtained in our main theorem is sharp. Understanding the sharp rates of decay of solutions is also necessary in general relativity for a proof of the strong cosmic censorship conjecture, see \cite{DLuk}.

\section{Notation and statement of the main theorem}\label{sec:notation}

\subsection{Notation} \label{subsec:implicit constants}
We write $X\ls Y$ to denote $|X| \leq CY$ for an implicit constant $C$ which may vary by line. Similarly, $X \ll Y$ will denote $|X| \le c Y$ for a sufficiently small constant $c>0$. In Sections 3 and 4, all implicit constants are allowed to depend only on the coefficients of $P$. In Sections 5-7 the constants may also depend on the initial data $\phi[0]$.

Next, we define notation for vector fields and some preliminaries. If $x =(x^1,x^2,x^3)\in\R^3$, we write 
\begin{align*}
r &:= \lr{ \sum_{i=1}^3 (x^i)^2 }^{1/2}, \quad u:= t-r,  \quad v := t+r.
\end{align*}
We write $\bo := -\pat^2+\Delta$. %

Let $\la a \ra$ denote a smooth function of $a$ equal to $1/2$ when $a\leq 1/2$ and $a$ if $a\geq 1$. This choice will be convenient later on when partitioning the space-time into finer regions.

In $\R^{1+3}$, we consider the three sets 
\[
\partial := \{\partial_t, \partial_1,\pa_2,\pa_3\}, \qquad \Omega := \{x^i \partial_j -
x^j \partial_i, \quad 1\leq i<j\leq 3\}, \qquad S := \{t \partial_t + \sum_{i=1}^3x^i \partial_{i}\},
\]
which are, respectively, the generators of translations, rotations and scaling. 
We now define the set of vector fields that we will use as commutators
\[Z := \{\pa,\Omega,S\} \] 
For a triplet $J=(i, j, k)$ we define $|J| = i + j + k$, $\phi_J$ is defined as $\phi_J(t, x) := (Z^J \phi)(t, x)$. More precisely,
\begin{equation}\label{vf defn}
\p_{J}(t,x) := (Z^J\phi)(t,x) := (\partial^i \Omega^j S^k \phi)(t,x), 
\end{equation}
$$\p_{\leq m} := \{\p_{J}\}_{|J|\le m}$$

We define the function class $$S^Z(f)$$ to be the collection of functions $w : \R \times \rt \to \R$ such that
\begin{equation}\label{SZdef} |Z^J w(t,x)| \ls_J |f| \end{equation}
whenever $J$ is a multiindex. We will frequently use $f = \jr^k$ for some real $k$. We also define $$S^Z_\text{radial}(f) := \{ w \in S^Z(f) : w \text{ is spherically symmetric} \}.$$

We also define 
$$S^Z_\text{der}(f) := \{ w \in S^Z(f): \jr \pa w \in S^Z(f)\}$$

We  denote by $\ang_i$ the angular derivatives
\begin{equation}
    \ang_i = \partial_i - \frac{x_i}{r} \partial_r,
\label{def:angular}
\end{equation}
and we define the set of tangential derivatives as
\[
\overline{\pa} := \{\pa_t + \pa_r, \ang_1, \ang_2, \ang_3\}.
\]
Notice that for any $C^1$ function $f : \R^3_x \to \R$,
$$\sum_{i=1}^3 |\ang_i f|^2 = |\nabla_x f|^2 - |\partial_r f|^2.$$

Tangential derivatives decay better near the light cone, as they satisfy
\begin{equation}\label{tangdef}
 \overline{\partial} \phi \in S^Z \Bigl(\frac{\ju}{r}\Bigr)\pa \p + S^Z \Bigl(\frac{1}{r}\Bigr) Z\p.
 \end{equation}
 
 Indeed, \cref*{tangdef} is a consequence of
 \[
 (\pa_t + \pa_r)\p = \frac1r S\p + \frac{r-t}{r}\pa_t\p,
 \]
 and 
\begin{equation}\label{eq:angular_relatedto_rotation}
     \ang_i\p = -\sum_{j\neq i}\frac{x_j}{r^{2}}\Omega_{ij}\p\in S^Z\Bigl(\frac{1}{r}\Bigr) \Omega\p.
 \end{equation}

Given a finite family of functions $\mathcal{G}$, we will also use the notation
 \[
 f \in S^Z(\la r\ra^k) \mathcal{G}
 \]
 to mean that
 \[
 f = \sum h_i g_i, \quad h_i\in S^Z(\la r\ra^k), \quad g_i\in \mathcal{G}.
 \]

For the initial data, we define the vector fields $Z|_{t=0}$ as $\{\pa_x, S|_{t=0} = ~r\pa_r, \Omega\}$. Given a function $f_0:\R^3\to\R$, we then define
\[
(f_0)_J = (Z|_{t=0})^J f_0.
\]

\subsection{Local energy norms}

We consider a covering $\R^3 = \bigcup_R A_R$ with the sets $A_R= \{R/2\leq \la r \ra \leq R\}$ for $R \geq 1$. %

Define the local energy norm $LE$ 
\begin{equation}
\begin{split}
 \| \p\|_{LE} &= \sup_{R \text{ dyadic}}  \| \la r\ra^{-\frac12} \p\|_{L^2 ([0,\infty) \times A_R)}  \\
 \| \p\|_{LE[t_0, t_1]} &= \sup_{R \text{ dyadic}}  \| \la r\ra^{-\frac12} \p\|_{L^2 ([t_0, t_1] \times A_R)},
\end{split} 
\label{ledef}\end{equation}
its $H^1$ counterpart
\begin{equation}
\begin{split}
  \| \p\|_{LE^1} &= \| \nabla_{t,x} \p\|_{LE} + \| \la r\ra^{-1} \p\|_{LE} \\
 \| \p\|_{LE^1[t_0, t_1]} &= \| \nabla_{t,x} \p\|_{LE[t_0, t_1]} + \| \la r\ra^{-1} \p\|_{LE[t_0, t_1]},
\end{split}
\end{equation}
as well as the dual norm
\begin{equation}
\begin{split}
 \| f\|_{LE^*} &= \sum_{R \text{ dyadic}}  \| \la r\ra^{\frac12} f\|_{L^2 ([0,\infty)\times A_R)} \\
 \| f\|_{LE^*[t_0, t_1]} &= \sum_{R \text{ dyadic}}  \| \la r\ra^{\frac12} f\|_{L^2 ([t_0, t_1] \times A_R)}.
\end{split} 
\label{lesdef}\end{equation}

Then we have the following scale invariant local energy estimate on Minkowski backgrounds:
\begin{equation}\label{localenergyflat}
\|\nabla_{t,x} \p\|_{L^{\infty}_t L^2_x} + \| \p\|_{LE^1}
 \lesssim \|\nabla_{t,x} \p(0)\|_{L^2} + \|\Box \p\|_{LE^*+L^1_t L^2_x}
\end{equation}
and a similar estimate involving the $LE^1[t_0, t_1]$ and $LE^*[t_0, t_1]$ norms.

The first estimate of this kind was obtained by Morawetz for the Klein-Gordon equation \cite{M}. There are many similar results obtained in the case of
small perturbations of the Minkowski space-time; see, for example, \cite{KSS}, \cite{KPV}, \cite{SmSo}, \cite{St}, \cite{Strauss},
\cite{Al2}, \cite{MS} and \cite{MT}. Even for large perturbations, in the absence of trapping, \cref*{localenergyflat} still sometimes holds, see for instance \cite{BH}, \cite{MST}.
In the presence of trapping, \cref*{localenergyflat} is known to fail, see \cite{Ral}, \cite{Sb}. 
We will assume that a similar estimate holds, with no loss of derivatives, for our operator $P$ after commuting with vector fields in $Z$.

\begin{definition}[Energy bounds and Local energy decay] \label{LED}
We say that $P$ has the strong local energy decay (SLED)
property if the following estimate holds for all $m \geq 0$, and $0\leq T_0<T_1 \leq \infty$:
\begin{equation}\begin{split}
\|\pa\pm\|_{L^\infty_t L^2_x([T_0,T_1)\times\R^3)}+  \| \p_{\le m}\|_{LE^{1}([T_0,T_1)\times\rt)} 
\\ \ls_m  \|\pa \p_{\le m}(T_0)\|_{\lt(\rt)} + \|(P\p)_{\le m}\|_{(\lolt + LE^*)([T_0,T_1)\times\rt)}.
\end{split}\label{eq:LED}\end{equation}
Here the implicit constant may depend on $m$, but not $T_0$ and $T_1$.
\end{definition}

For example, \cref*{LED} holds if the operator $P$ is a small perturbation of $\Box$, see for instance \cite{MT}. More generally, strong local energy estimates (for $m=0$) were shown in \cite{MST}, provided there are no negative eigenvalues or real resonances. Provided that $P$ is stationary, one can extend the result to \cref*{eq:LED} by commuting with the vectors fields in $Z$. 

\subsection{Statement of the main theorem}

Let $h=g-m$, where $m$ denotes the Minkowski metric. Let $\xs \in (0,\infty)$ be real. We make the following assumptions on the coefficients of $P$: 
\begin{equation}\label{coeff.assu}
\begin{split}
 h^{\alpha\beta} \in S^Z(\jr^{-1-\xs}), \quad  V\in S^Z(\jr^{-2-\xs})\\
B^{\alpha} \in S^Z_{der}(\jr^{-1-\xs})   \\
g^\xo \in S^Z_\text{radial}(\jr^{-2-\xs})
\end{split}
\end{equation}

We assume that $Q$ is a generalized null form, in the sense that 
\begin{equation}\label{Qass}
Q[\pa\phi, \pa\phi] \in S^Z(1) \pa\phi \tpa\phi.
\end{equation}

\begin{theorem}[Main theorem]\label{thm:main} Assume that $P$ has the SLED property (\cref*{LED}), and that the coefficients of $P$ satisfy \cref*{coeff.assu}. Consider the initial-value problem
\begin{equation}\label{eq:problem}
P\phi = Q[\pa\phi, \pa\phi], \qquad (\p(0,\cdot), \pat\p(0,\cdot)) = (\phi_0,\phi_1)
\end{equation}
where $Q$ is a generalized null form.
\begin{enumerate}[label=(\roman*)]
\item
Assume that $(\phi_0,\phi_1) \in H^{15}(\rt) \times H^{14}(\rt)$. Then there is a constant $\eps_0> 0$ so that, if
$$ \sum_{|J|=0}^{14}\|\pa_x(\phi_0)_J (\cdot)\|_{\lt(\R^3)} + \|(\phi_1)_J(\cdot)\|_{L^2(\R^3)} \leq \eps_0,$$
then \cref*{eq:problem} has a unique global solution forward in time. 
\item
 Let $m \in\N$. Then there exists an integer $N>m$ such that, additionally assuming there is an absolute constant $C > 0$  with the initial data satisfying 
\begin{equation}\label{eq:maintheoreminitialdata}
    \sum_{|J|=0}^N \|(1+|\cdot|)^{1/2}( \phi_0)_J(\cdot)\|_{L^2(\R^3)} + \|(1+|\cdot|)^{3/2}(\phi_1)_J(\cdot)\|_{L^2(\R^3)}  \leq C,
\end{equation}
the unique forward solution and its vector fields up to order $m$ satisfy
\begin{equation}\label{opt}
\sum_{|J|=0}^m |\phi_{J}(t,x)| \lesssim C \langle v \rangle^{-1} \langle u \rangle^{-1}.
\end{equation}
Moreover, 
\begin{equation}\label{eq:derivatives}
\sum_{|J|=0}^m |\partial \phi_{J}(t,x)| \lesssim C \langle r \rangle^{-1} \langle u \rangle^{-2}.
\end{equation}
\end{enumerate}
\end{theorem}

\begin{remark}[Sharpness of the decay rate]
The result \cref*{opt} is optimal for $\phi$ in our context for compactly supported smooth initial data. Heuristically, this comes from the fact that in general $Q$ cannot decay faster than $\frac{1}{r^3}$ near the light cone $r=t$. Indeed, even for solutions to the linear equation $\Box \tilde\phi = 0$, the best one can prove is that $Q[\pa\tilde\phi, \pa\tilde\phi] \lesssim \frac{1}{r^3}$ near the light cone. Solving the linear wave equation $\Box w = F$ with $F\approx r^{-3}$ near the cone only yields $t^{-2}$ decay in a compact region, which is precisely \cref*{opt}. 

To see this more concretely, consider a simple model on Minkowski space that falls within our general framework:  \begin{equation*} \Box \phi = c(t, x) ((\partial_t\phi)^2 - |\partial_x\phi|^2). \end{equation*} 
Fix $N, M \gg 1$. The coefficient $c(t,x)$ is a smooth function such that
\begin{equation*}
c(t, x) = 
\begin{cases} 
M - \langle t \rangle^{-N} & \text{for } |x| < 1, \\
M - \langle t-r \rangle^{-N} & \text{for } |x| > 2.
\end{cases}
\end{equation*}
with a smooth interpolation on $1 \le |x| \le 2$. 
\end{remark}

\subsection{Space-time partition}\label{subsec:spacetime_partition}
We decompose the set $\R\times \R^3$ based on either distance from the cone $\{r=t\}$ or distance from the origin $\{r=0\}$. Given a number $T>8$, let $I_T = [\frac34 T, T]$, and $C_T :=  I_T\times \R^3$. We partition $C_T$ into dyadic regions $C^R_T$ and $C^U_T$, depending on whether we are far or close to the light cone $\{r=t\}$, by starting on a portion of $\{(t,x) : t=\frac 34 T\}$, and flowing along the vector field $\frac1t S$. More precisely, observe that
\[
\gamma_{(t,x)}(s) = (t+s, \frac{t+s}{t} x),
\]
the straight line connecting the origin to $(t,x)$, is the integral curve of $\frac1t S$ starting at $(t,x)$. 

We will define two types of trapezoidal regions that respect the flow of $\frac1t S$. Far away from the cone, let
\begin{equation}
\label{eq:setsR}\begin{split}
& C^R_T := \{\gamma_{(\frac34 T,x)}(s): R/2\leq |x| \leq R, 0\leq s\leq \frac 14 T\} = C_T \cap \Bigl\{\frac{R/2}{\frac34 T} \leq  \frac{ r}{t} \leq \frac{R}{\frac 34 T}\Bigr\} , 
\\ & \quad  2 \leq R \, \,\text{and } \, \left( R \leq \frac{T}{2} \, \text{or} \, R\geq 3T \right). 
\\ & C^R_T := \{\gamma_{(\frac34 T,x)}(s): 0 \leq |x| \leq 1, 0\leq s\leq \frac14 T\} = C_T \cap \Bigl\{0 \leq \frac{r}{t} \leq \frac{1}{\frac 34 T}\Bigr\}, \quad R=1.
\end{split}\end{equation}

Near the cone, we define the following sets. 

\begin{equation}\label{eq:setsU}
\begin{split}
& C^U_T := \{\gamma_{(\frac34 T,x)}(s): \frac34 T-U \leq |x| \leq \frac34 T-\frac12 U, 0\leq s\leq \frac14 T\} = C_T \cap \Bigl\{\frac{\frac34 T-U}{\frac34 T} \leq \frac{r}{t}  \leq \frac{\frac34 T-\frac12 U}{\frac 34 T}\Bigr\}, \\& 2\leq U \leq \frac{T}{4} \\
& C^U_T := \{\gamma_{(\frac34 T,x)}(s): \frac34 T-1 \leq |x| \leq \frac34 T+1, 0\leq s\leq \frac14 T\} = C_T \cap \Bigl\{\frac{\frac34 T-1}{\frac34 T} \leq \frac{r}{t} \leq \frac{\frac34 T+1}{\frac 34 T}\Bigr\}, \\& \quad U=1 \\
& C^U_T := \{\gamma_{(\frac34 T,x)}(s): \frac34 T-\frac12 U \leq |x| \leq \frac34 T-U, 0\leq s\leq \frac14 T\} = C_T \cap\Bigl\{\frac{\frac34 T-U}{\frac34 T} \leq \frac{r}{t} \leq \frac{\frac34 T-2U}{\frac34 T}\Bigr\}, \\ & -2 \geq U \geq -\frac{3T}{4}
\end{split}
\end{equation}

We also define
\[
C_T^{\mathrm{cone}} = \bigcup_{\substack{U \text{ dyadic} \\ -\frac{3T}{4} \le U \le \frac{T}{4}}} C^U_T.
\]

Note that in the first set with $2 \leq U \le T/4$, the region $C^U_T$ corresponds to the set of points reached by flowing along the integral curves of $(1/t)S$ for time $T/4$, starting from points $(3T/4, x)$ where $3T/4 - U \le |x| \le 3T/4 - \frac12 U$; this is an area \textit{inside} the light cone. Note that the $U = 1$ case is special because it's the only one that includes the light cone $r = t$ at the initial time $3T/4$. Finally, part of $C^U_T$ for negative $U$ covers some areas \textit{outside} the light cone.

There are three reasons for the choice of these particular regions. The most important reason is that proving Sobolev embeddings in the regions $C^R_T$ and $C^U_T$ in Lemma~\ref{SobEmbExt} is most natural when looking at regions defined by integral curves of $S$. We also have that $C^R_T$ and $C^U_T$ have disjoint interiors (for $R$ and $U$ dyadic), and they cover $C_T$, i.e.
\[
 C_T = \bigcup_R \crt \cup \bigcup_U\cut
 \]

Finally, due to our numerology, we can enlarge the regions $C_T^R$ defined in \cref*{eq:setsR} twice to a rectangular region $\tilde C_T^R$; this will prove useful in some of the results in Section 5, in particular Lemma~\ref{L2improv}.  More precisely, let 
\begin{equation}\label{DTR}
\begin{split}
 D_T^R := C_T\cap\{\frac12 R\leq \la r\ra\leq \frac53 R\} \\ \tilde C_T^R := \Bigl\{\frac34 T -\frac1{100}T \leq t\leq T+ \frac1{100}T, \frac12 R-\frac{1}{100} R\leq \la r\ra\leq \frac53 R+ \frac{1}{100}R\Bigr\}    
\end{split}\end{equation}
This enlarged region $\tilde C^R_T$ is used to give some ``breathing room" around $C_T^R$ for Lemma~\ref{L2improv}. The small extensions (by factors of $\frac1{100}$) ensure that $\tilde C_T^R$ properly contains $C_T^R$ without significantly altering its properties. 

We observe that we have 
\[
C^R_T \subset D_T^R \subset \tilde C^R_T
\]
and thus in $C^R_T$ (and also the larger region $\tilde C^R_T$):
\begin{equation}\label{eq:TUest}
|t-r|\sim t+r \sim \max\{R, T\}.
\end{equation}

Similarly, if we choose 
\begin{equation}\label{DTU}
\begin{split}
D^U_T := C_T \cap\{\frac12 |U|\leq \la t-r\ra\leq \frac53 |U|\} \\ \tilde C^U_T := \Bigl\{\frac34 T -\frac1{100}T \leq t\leq T+ \frac1{100}T, \quad \frac12 |U|-\frac{1}{100}|U|\leq \la t-r\ra\leq \frac53 |U|+ \frac{1}{100}|U|\Bigr\}  
\end{split}\end{equation}
we have that
\[
C^U_T \subset D_T^U \subset \tilde C^U_T
\]
and thus in $C^U_T$ (and also the larger region $\tilde C^U_T$):
\begin{equation}\label{TRest}
\la t-r\ra\sim |U|, \quad r\sim t\sim T.
\end{equation}

Solutions to wave equations on variable-coefficient space-times exhibit the property that their derivatives in general\footnote{A subcollection of the full set of derivatives decay more rapidly, but here we refer to generic derivatives.} decay differently depending on whether the solution is evaluated near $r=0$ or near $r=t$. Thus we use the $C^R_T$ and $C^U_T$ sets to distinguish between these situations respectively. This heuristic is also true for nonlinear wave equations on variable-coefficient space-times; see \cref*{derbound} and analogous pointwise decay estimates for derivatives of solutions to nonlinear wave equations in \cite{Looi2,Looi,Looi22Improved}.

\subsection{Exterior region} \label{subsec:ext}

In order to obtain optimal estimates in the exterior region, we need to also define regions on the same scale in space and time. We will \textit{not} use these regions in the proof of global existence. However, once global existence is proved, we can use them to improve decay in the region $r\gg t$.

Let $\tilde\gamma_{(t,x)}(s) = (\frac{t}{|x|}(|x|+s), \frac{x}{|x|}(|x|+s))$ denote the flow of $\frac1r S$ stating at $(t,x)$. We define for $T > 1$ and $R\geq \frac32 T$
\[\begin{split}
C_R^T := \{\tilde\gamma_{(t,x)}(s): |x|=R, T\leq t \leq T+R/10, 0\leq s\leq R\} \\ = \{R\leq r\leq 2R, \quad \frac{T}{R} \leq\frac{t}{r} \leq \frac{T+R/10}{R}\}
\end{split}\]

This region is similar to $C^R_T$, but ``with $r$ and $t$ having swapped roles.'' Thus it makes sense to flow from an $r=$constant slice (rather than $t=$constant slice) along the flow lines of $\frac1r S$ (rather than $\frac1t S$). If we perform Sobolev embedding in $C^T_R$, we obtain a scaling factor in the pointwise estimate $|\phi(t,x)| \lesssim \la x \ra^{-1/2} \|\phi_{\le k}\|_{LE^1} \approx R^{-1/2} \|\phi_{\le k}\|_{LE^1}$. That this is a power of $R$, $R^{-1/2}$, rather than a power of $T$, is important.

As with the previous regions, we will need to extend the regions $C_R^T$. We define
\begin{equation}\label{DRT}
\begin{split}
 D^T_R := \{T\leq t\leq T+R/10, \quad R\leq r \leq 2R\} \\ \tilde C_R^T := \Bigl\{T -\frac1{10}R \leq t\leq T+ \frac1{10}R, \quad R-\frac{R}{10}\leq r \leq 2R+ \frac{R}{10}\Bigr\} 
\end{split}\end{equation}

Once again, we have
\[
C_R^T\subset D_R^T\subset \tilde C_R^T
\]
and in $\tilde C_R^T$ 
\begin{equation}\label{RTest}
 r\sim t+r\sim \la t-r\ra\sim R.
\end{equation}

Note however that we do not have $t\sim T$ in the extended region $\tilde C_R^T$. when $R\gg T$. Actually, the region ${\tilde C}_R^T$ extends to $t<0$ when $R\gg T$. This will not be an issue, since we will only use this region to establish Lemma~\ref{L2improv}, and restrict it to $t\geq 0$.

\subsection{Interior region} \label{subsec:int}

In order to obtain estimates in the interior region in Section 5, we define
 \[
\inte := \bigcup_{\substack{R \text{ dyadic} \\ R \le T/2 \text{ or } R=1}} \crt = \left\{ (t,y): \frac34T\leq t\leq T,\quad  |y| \leq \frac{2t}3 \right\} %
\]
We also define its enlargement (this time only in space) by
\[
\tinte := \left\{ (t,y)= \gamma_{(\frac34T,x)}(s), \quad 0\leq s\leq \frac14T, \quad |x| \leq\frac12T \right\} = \left\{ \frac34T\leq t\leq T,\quad  |y| \leq \frac{3t}4 \right\} 
\]

\subsection{Notation for the symbols $n$ and $N$} \label{subsec:N}
Throughout the paper the integer $N$ will denote a fixed and sufficiently large positive number, signifying the highest total number of vector fields that will ever be applied to the solution $\p$ to \cref*{eq:problem} in the paper.

We use the convention that the value of $n$ may vary from line to line.%

\section{Pointwise bounds from local energy}\label{sec:pointwise}

In this section we will show that local energy bounds imply certain weak pointwise bounds, see Proposition~\ref{inptdcExt} and Proposition~\ref{DerProp}. These bounds are sufficient to prove global existence in Section 4.

We start with the following Klainerman-Sideris type estimate for the second derivative.
\begin{lemma}\label{2ndDeBd'}
Assume $\p$ is sufficiently regular. We then have for all $r\gg1$
 \begin{equation}\label{2ndD}
|\pa^2\p_J| \ls
\left(\f1\jr+\f1\ju\right)|\pa\p_{\le|J|+1}| +  \left( 1 + \f{t}\ju \right) \jr^{-2}|\p_{\le|J|+2}|
+ \left( 1 + \f{t}\ju \right) |(P\p)_{\le|J|}|.
 \end{equation}
\end{lemma}

\begin{proof}

Note first that
\begin{equation}\label{2nd'sEst}
|\pa^2\p_J| \ls
\left(\f1\jr+\f1\ju\right)|\pa\p_{\le|J|+1}|
+ \left( 1 + \f{t}\ju \right) |(\Box\p)_{\le|J|}|.
 \end{equation}

The case $|J|=0$ is an immediate consequence of Lemma 2.3 from \cite{KS}. The general case follows after commuting with vector fields.

It is thus enough to estimate the difference $P-\Box$. We write by \cref*{P def}
\[
P - \Box = h^{\alpha\beta}\pa_{\x\xb} + (\paa\hab)\pab + g^\xo(t,x) \Delta_\xo + B^\x(t,x)\paa + V(t,x)
\]
Using the assumptions on the coefficients in \cref*{coeff.assu} we have
\begin{align*}
(P - \Box)\p \in S^Z(\jr^{-1-\xs})(\pa^2\p + \pa\p)
+ S^Z(\jr^{-2-\xs})\Omega^{\leq 2}\p  
\end{align*} 

After applying vector fields we thus obtain
\begin{equation}\label{error}
\left|\left((P - \Box)\p\right)_{\leq |J|}\right| \in S^Z(\jr^{-1-\xs})|\pa^2\p_{\leq |J|}| + S^Z(\jr^{-1-\xs}) |\pa\p_{\leq |J|}| + S^Z(\jr^{-2-\xs}) |\p_{\leq |J|+2}|
\end{equation}

The conclusion now follows from \cref*{2nd'sEst} and \cref*{error}, since  the first term on the RHS of \cref*{error} can be absorbed in the LHS of \cref*{2nd'sEst} for $r \gg 1$. %
\end{proof}

The main tool for turning local energy estimates into pointwise bounds is the following lemma. For the proof of global existence in Section~\ref{sec:existence}, it will be important that the norms on the right-hand sides in Lemma~\ref{SobEmbExt} do not have enlarged domains (in time). While a proof of this lemma with enlarged domains on the right-hand side is possible, we seek a proof without an enlarged domain on the right-hand side, which leads us to use the idea of defining the dyadic regions using the flow of vector fields (see subsection~\ref{subsec:spacetime_partition} for the definitions).
\begin{lemma}[Dyadically localized bounds]
\label{SobEmbExt}
Let $w\in C^4$. We then have 
 \begin{equation}\label{U-Sob}
 \| w\|_{L^\infty(\cut)} \ls \sum_{i\le 1,j\le 2} \f1{(T^3|U|)^{1/2}}  \|\Omega^j S^i  w\|_{L^2( C^U_T)} +  \lr{\f{|U|}{T^3}}^\f12 \|\pa \Omega^j S^i  w\|_{L^2( C^U_T)}.
\end{equation}

\begin{equation}\label{R-Sob}
 \begin{split}   
& \| w\|_{L^\infty(C^R_T)} 
  \ls \sum_{i\le 1,j\le 2}
\f1{(R^3T)^{1/2}}
    \|\Omega^j S^i  w\|_{L^2(C^R_T)} +  
    \f1{(RT)^{1/2}} \|\pa \xO^j S^i   w\|_{L^2( C^R_T)}, \quad R\geq 2. 
\\& 
 \| w\|_{L^\infty(C^R_T)} 
  \ls \sum_{i\le 1,j\le 2}
\f1{T^{1/2}}
    \|\pa^j S^i  w\|_{L^2(D^R_T)}, \quad R=1.   
\end{split}    
\end{equation}

\begin{equation}\label{R-Sob2}
    \| w\|_{L^\infty(C^T_R)} 
  \ls \sum_{i\le 1,j\le 2}
\f1{R^2}
    \|\Omega^j S^i  w\|_{L^2(C^T_R)} +  
    \f1{R} \|\pa \xO^j S^i   w\|_{L^2( C^T_R)}. 
\end{equation}

\begin{equation}\label{cone-Sob}
 \| w\|_{L^\infty(C_T^{\mathrm{cone}})} \ls \sum_{i\le 1,j\le 2} \f1{T^2}  \|\Omega^j S^i  w\|_{L^2( C_T^{\mathrm{cone}})} +  {\f{1}T} \|\pa \Omega^j S^i  w\|_{L^2( C_T^{\mathrm{cone}})}.
\end{equation}
\end{lemma}

\begin{proof}
Let us consider the proof for $C_T^R$. For any fixed point $(t, x=r\omega) \in C_T^R$, define the curve:

\begin{equation}\label{eq:gamma}
\gamma_{T,x}(s) = (\tau(s), \rho(s)\omega) = \left(T+s, \frac{r}{t}(T+s)\omega\right)
\end{equation}
with $\tau(s) = T+s, \rho(s) = \frac{r}t (T+s)$. Thus $\tau(0) = T$ and $\dot \tau(s) = 1$.

Let $\gamma_{T,x}(s) = (T+s, \frac{r}{t}(T+s)\omega)$ be the integral curve of $\frac{1}{\tau(s)}S$ passing through $(t,r\omega)$ at $s=t-T$. We use the following key inequality, which holds for any interval $I$, any function $f \in C^1(I)$ and any point $x \in I \subset \mathbb{R}$. It is a one-dimensional Poincaré inequality.
\begin{equation}\label{eq:1dest}
f^2(x) \lesssim \frac{1}{|I|}\int_I |f(y)|^2 \ dy + |I|\int_I |f'(y)|^2  \ dy
\end{equation}

Applying \cref*{eq:1dest} along the integral curves of $S$ and using the fact that $\frac{d}{ds} w(\gamma_{T,x}(s)) = \frac{1}{T+s} Sw(\gamma_{T,x}(s))$, we obtain:
\begin{equation}\label{eq:w_bound}
|w(t, r\omega)|^2 \lesssim T^{-1}\int_{-\frac{T}4}^{0} |w(\gamma_{T,x}(s))|^2 + |Sw(\gamma_{T,x}(s))|^2 ds
\end{equation}

Next, we integrate in the $r$-direction. Define 
\[
I_s = C_T^R \cap\{\tau= T+s\}.
\]
More precisely,
$$I_s = \Bigl\{\rho\in\R : \frac{R/2}{\frac34 T} \leq \frac{\rho}{T+s} \leq \frac{R}{{\frac34 T}}\Bigr\}, \quad -\frac{T}{4} \leq s \leq 0, \quad R\geq 2 $$ 
$$I_s = \{\rho\in\R : 0 \leq \frac{\rho}{T+s} \leq \frac{1}{\frac34 T}\Bigr\}\}, \quad -\frac{T}{4} \leq s \leq 0, \quad R=1. $$

The case $R=1$ follows now simply from observing that 
$$I_s\subset \{|x|\leq \frac43\} \subset\subset \{|x|\leq \frac83\}$$
and applying Sobolev embeddings to $w$ and $Sw$ on $\{|x|\leq \frac43\}$.

In the case $R\geq 2$, note first that  $|I_s| \sim R$ and $\rho \sim R$ in $I_s$. Applying \cref*{eq:1dest} to $w(T+s, \rho\omega)$ and $Sw(T+s, \rho\omega)$ as functions of $\rho$ on the interval $I_s$:
\begin{equation} 
\begin{split}
|w(\gamma_{T,x}(s))|^2 \lesssim R^{-1}\int_{I_s} |w(T+s,\rho\omega)|^2 d\rho + R \int_{I_s}|\partial_\rho w(T+s, \rho\omega)|^2 d\rho \\
|Sw(\gamma_{T,x}(s))|^2 
\lesssim R^{-1}\int_{I_s} |Sw(T+s, \rho\omega)|^2 d\rho + R \int_{I_s}|\partial_\rho Sw(T+s, \rho\omega)|^2 d\rho
\end{split}
\end{equation}

Since $\rho \sim R$ in $I_s$, we can replace $d\rho$ with $R^{-2}\rho^2 d\rho$:

\begin{equation}\label{wS} 
\begin{split}
|w(\gamma_{T,x}(s))|^2 &\lesssim R^{-3} \int_{I_s} |w(T+s, \rho\omega)|^2 \rho^2 d\rho + R^{-1} \int_{I_s}|\partial_\rho w(T+s, \rho\omega)|^2 \rho^2 d\rho \\
|Sw(\gamma_{T,x}(s))|^2 &\lesssim R^{-3} \int_{I_s} |Sw(T+s, \rho\omega)|^2 \rho^2 d\rho + R^{-1} \int_{I_s}|\partial_\rho Sw(T+s, \rho\omega)|^2 \rho^2 d\rho
\end{split}
\end{equation}

Finally, we apply Sobolev embedding on the sphere 
\[
|f(\rho\omega)|^2 \lesssim \sum_{j=0}^2 \|\Omega^j f(\rho\cdot)\|^2_{L^2(\S^2)}
\]
to $w$, $Sw$, $\pa_\rho w$ and $\pa_\rho Sw$ in \cref*{wS}.

Substituting into \cref*{eq:w_bound} and changing the order of integration, noting that $\bigcup_{s \in [-\frac{T}{4},0]} I_s \times \S^2= C_T^R$, we obtain  \cref*{R-Sob}.

The proofs of \cref*{U-Sob}  
 and \cref*{cone-Sob} are almost identical, except that the integral in the $r$-direction is on an interval of size comparable to $|U|$ or $T$ respectively. Finally, in the proof of \cref*{R-Sob2} we integrate in the $t$-direction and in the $\frac1r S$ direction on intervals of size comparable to $R$.

For concreteness, we sketch the ideas of the proof of \cref*{cone-Sob} and highlight what is different in it.   We apply \cref*{eq:1dest}
first along integral curves of $S/t$,  and second along the radial direction using spatial scale $L=T$ in both cases. Using the fact that $r \sim T$ in $C_T^{\text{cone}}$, we obtain 
\begin{equation} 
\begin{split}
& |w(t, r\omega)|^2 \lesssim T^{-1}\int_{-\frac{T}4}^{0} |w(\gamma_{T,x}(s))|^2 + |Sw(\gamma_{T,x}(s))|^2 ds  \\ &\lesssim T^{-2} \int_{-\frac{T}4}^{0}\int_{I_s} \Bigl(|S^{\leq 1} w(T+s, \rho\omega)|^2  +  |\partial_\rho S^{\leq 1} w(T+s, \rho\omega)|^2 \Bigr)\rho^2 d\rho 
\end{split}
\end{equation}

 Combining with Sobolev on $\mathbb{S}^2$ and spacetime volume factors ($T^{-2}$) leads to $$\|w\|_{L^\infty}^2 \lesssim T^{-4}\|(I, S)W\|_{L^2}^2 + T^{-2}\|\partial (I, S)W\|_{L^2}^2,$$ giving \cref*{cone-Sob} after taking the square root.
\end{proof}

The following lemma, which is similar to Lemma 9.1 in \cite{LR2}, will be applied ``near the cone'' $|x|=t$.

\begin{lemma}\label{Hardy}
Let $A[a,b] := \{ x \in \R^3 : a \le |x|\le b\}.$
If $f\in C^1(\R \times \R^3;\R)$, $t \geq 1$, and $\delta \in(0,1)$, then
\begin{align}\label{con.Hardy}
\begin{split}
\int_{A[t/2,3t/2]} \nm^{-1+\delta} f(t,x)^2 dx &\ls t^{2\delta}\int_{A[t/4,7t/4]}\la t-r\ra^{1-\delta} |\pa_r f(t,x)|^2 dx \\
&\ + \f1{t^{1-\delta}} \left(\int_{A[t/4,t/2]} f(t,x)^2 dx + \int_{A[3t/2,7t/4]} f(t,x)^2 dx\right)
\end{split}
\end{align}
\end{lemma}

\begin{proof}Fix $t\geq 1$ and $\omega\in\S^{2}$, and set
\[
g(r):=f\bigl(t,r\omega\bigr), \qquad 0<r<\infty .
\]
Because $dx=r^{2}\,dr\,d\omega$, the desired estimate will follow once we
prove the corresponding inequality for $u$ (with the same constants) and
then integrate over $\S^{2}$.

Let $\chi : [0,\iy) \to [0,1]$ be a cutoff such that $\chi(s) = 1$ for $1/2 \le s \le 3/2$ and 0 when $s \le 1/4$ and $s \ge 7/4$. We will show that
\begin{equation}\label{P1d}
\int\nm^{-1+\delta} \chi^2(r/t) g(r)^2 r^2dr \ls  t^{2\delta}\int \nm^{-1+\delta} |g'(r)\chi(r/t)|^2 r^2dr  + \f1{t^{1-\delta}} \int |g(r) \chi'(r/t)|^2 r^2 dr.
\end{equation}
 
We have
$$g(r)^2\chi^2(r/t) - g(7t/4)^2\chi^2((7t/4) / t) = -2\int_r^{7t/4} g(\rho) \chi(\rho/t) \cdot (g(\rho) \chi(\rho/t))'d\rho.$$
Recall that $\chi(7t/4) = 0$. Hence 
$$g(r)^2 \chi^2(r/t)  \ls  \int_r^{7t/4} |g(\rho) \chi(\rho/t) \cdot  ( g(\rho) \chi(\rho/t))' | d\rho$$
 
 Let $w(r)=\la t-r\ra^{-1+\delta}$ and $H(r) = g(r)\chi(r/t)$. Multiplying the inequality by $w(r)r^2$, integrating in $r$ from $t/4$ to $7t/4$ and applying Fubini's theorem yields:
\begin{align}
\int_{t/4}^{7t/4} w(r) r^2 H^2(r) dr &\ls  \int_{t/4}^{7t/4} w(r) r^2 \left( \int_r^{7t/4} |H(\rho)| |H'(\rho)| d\rho \right) dr \nonumber \\
&=  \int_{t/4}^{7t/4} |H(\rho)| |H'(\rho)| \left( \int_{t/4}^{\rho} w(r) r^2 dr \right) d\rho, \label{eq:after_forder}
\end{align}

Since $r\approx \rho \approx t$, and, crucially, 
\[
\int_{t/4}^{\rho} w(r) dr \lesssim t^{\delta}
\]
we obtain
\begin{align*}
\int_{t/4}^{7t/4} w(r) H^2(r) r^2dr &\ls \int_{t/4}^{7t/4}  t^{\delta} |H(r)| | H'(r)| r^2dr
\end{align*}%
Using Cauchy-Schwarz, we obtain
\begin{equation}\label{H'}
\int_{t/4}^{7t/4} w(r) H^2(r) r^2dr \ls \int_{t/4}^{7t/4} t^{2\delta} w(r)^{-1} |H'(r)|^2 r^2 dr.
\end{equation}
By the chain rule, $\pa_r(\chi(r/t)) \ls |\chi'(r/t)| \cdot \f1t$, and thus
\[
|H'(r)| \ls |g'(r)\chi(r/t)| + t^{-1}|g(r) \chi'(r/t)|
\]

Plugging the above into \cref*{H'} yields
\begin{align*}
\int_{t/2}^{3t/2}\nm^{-1+\delta} H^2(r) r^2dr
&\ls  \int_{t/4}^{7t/4}t^{2\delta}\nm^{1-\delta}|g'(r)\chi(r/t)|^2 r^2dr \\& + \f1{t^{1-\delta}} \int_{t/4}^{7t/4} |g(r) \chi'(r/t)|^2 r^2 dr,
\end{align*}
which is \cref*{P1d}.
\end{proof}

Our next proposition yields global pointwise bounds for $\phi_J$ under the assumption that the local energy norms are finite. These estimates are sharp from that point of view, but can be improved for solutions to \cref*{eq:problem}, see Sections 5-7.

\begin{proposition} \label{inptdcExt}
Let $T\geq 8$ be fixed and $\p$ be any sufficiently regular function.
There is a fixed positive integer $k$, %
such that for any multi-index $J$ with $|J|\le N - 3$, we have:
\begin{align}\label{u/v decay}
\begin{split}
|\p_{J}(t, x)| \leq \bar C_{|J|} \|\p_{\leq |J|+3}\|_{LE^1(I_T)} \la t\ra^{-1/2}, \quad\forall (t,x)\in C_T. %
\end{split}\end{align}
\end{proposition}

\begin{proof}
Away from the cone, meaning in the $C_T^R$ regions, the bound \cref*{u/v decay} is a straightforward consequence of \cref*{R-Sob}. 

In the region $C_T^{\mathrm{cone}}$ (i.e. ``near the cone''), \cref*{u/v decay} is an immediate consequence of \cref*{cone-Sob}, since $C_T^{\mathrm{cone}}$ lies inside a single (spatial) dyadic region.
\end{proof}

We now obtain an improved bound on the derivatives. Note that here it is crucial that $\phi$ is a solution to \cref*{eq:problem} , since we need to use \cref*{2ndDeBd'}.
Let 
\begin{equation}
    \label{def:mu}
    \mu := \mu(t,r) := \min(\jt,\ju)^{1/2}.
\end{equation}
\begin{proposition}\label{DerProp}
Let $T$ be fixed and $\p$ solve \cref*{eq:problem} in $C_T$. Then for any region $C\in\{C_T^R, \cut\}$ and $m\geq 0$ we have %
\begin{equation}\label{deBound}
 \|\pa\pm\|_{L^\iy(C)} \le \bar C_{m} \f1\mu \left(\f1{\jr} + \|\pa\p_{\le \lceil \frac{m+3}2\rceil }  \|_{L^{\infty}(C)}\right)\|\p_{\le m+5} \|_{LE^1(I_T)} 
\end{equation}
\end{proposition}

\begin{proof} 
Note first that if $r\lesssim 1$, the bound follows immediately from \cref*{SobEmbExt}. Subsequently, we assume that $r\gg1$. 

Note first that
\[ 
(P\phi)_{\leq m+3} \ls \left| \left(\pa\p_{\leq \lceil \f{m+3}2 \rceil}\right)\left(\pa\p_{\le m+3} \right) \right|
\]

We also have
\begin{enumerate}
    \item In $C^R_T$, the bound in \cref*{2ndD} is 
\begin{equation}\label{2ndDfc}
\pa^2\p_{\leq m+3} \ls  R\inv|\pa\p_{\leq m+4}| + R^{-2}|\p_{\leq m+5}|+|(P\p)_{\le m+3}|.
\end{equation}
\item In $C^U_T$, the bound in \cref*{2ndD} is 
\begin{equation}\label{2ndDnc}
\pa^2\p_{\leq m+3} \ls |U|\inv|\pa \p_{\leq m+4}| +  \f1{T|U|} |\p_{\leq m+5}|+ T|U|\inv |(P\p)_{\le m+3}|.
\end{equation}
\end{enumerate}

We now apply \cref*{SobEmbExt} in our region $C$. When $C=C_T^R$ we obtain, using \cref*{R-Sob} and \cref*{2ndDfc}:
\[\begin{split}
 \|\pa\p_{\le m} \|_{L^{\infty}(C_T^R)} & \lesssim \frac{1}{(R^3T)^{1/2}}\|\pa\p_{\le m+3}\|_{L^2(C_T^R)} + \frac{1}{(RT)^{1/2}}\|\pa^2\p_{\le m+3}\|_{L^2(C_T^R)} \\ 
& \ls \frac{1}{RT^{1/2}} \|\p_{\le m+5} \|_{LE^1(I_T)} + \frac{\|\pa\p_{\leq \lceil\f{m+3}2 \rceil}\|_{L^{\infty}(C_T^R)}}{(RT)^{1/2}} \|\pa\phi_{\leq m+3}\|_{L^2( C^R_T)} \\
& \ls \left(\frac{1}{RT^{1/2}} + T^{-1/2}\|\pa\p_{\leq\lceil \f{m+3}2 \rceil}\|_{L^{\infty}(C_T^R)} \right) \|\p_{\le m+5} \|_{LE^1(I_T)}.
\end{split}\]
Here, $(RT)^{-1/2}\|\pa\p_{\leq \lceil\f{m+3}2 \rceil }\|_{L^{\infty}(C_T^R)} \|\pa\phi_{\leq m+3}\|_{L^2( C^R_T)}$ comes from the $(P\phi)_{\le m+3}$ term. We have placed the $\pa \phi_{\le m+4}$ term from \cref*{2ndDfc} into the $\|\phi_{\le m+5}\|_{LE^1}$ term.

When $C=C_T^U$ we similarly obtain by \cref*{U-Sob} and \cref*{2ndDnc}:
\[
 \|\pa\p_{\le m} \|_{L^{\infty}(C_T^U)} \ls \left(\frac{1}{T|U|^{1/2}} + |U|^{-1/2}\|\pa\p_{\leq \lceil\f{m+3}2 \rceil}\|_{L^{\infty}(C_T^U)} \right) \|\p_{\le m+5} \|_{LE^1(I_T)}
\]
This completes the proof.
\end{proof}

\section{The proof of small data global existence}\label{sec:existence}
In this section, we prove our first theorem. For any $N\in\N$, define
\[
\mathcal{E}_N(t) = \|\partial  \phi_{\leq N}\|_{L^{\infty}([0, t]; L^2(\R^3))} + \|\phi_{\leq N}\|_{LE^1([0, t]\times \R^3)}.
\]

We also define $\ti N = N-2$, and $N_1 = \lfloor N/2\rfloor$. 

The following is a black box theorem applicable to many circumstances. Its proof uses a detailed bookkeeping required to implement the bootstrap argument.

\begin{theorem}\label{thm:semi}
Let $N\ge 14$ be an integer, and fix $0<\delta \ll 1$ (say, $\delta = 1/100$). Then there exists $\eps_0 > 0$ so that, if the initial data satisfies
\[
\calE_N(0) \le \eps_0,
\]
then the solution to equation \cref*{eq:problem} exists globally in time. Moreover, there exist positive constants $\nu_N$ and $\tilde C$ so that the solution satisfies
\begin{equation}\label{GE1}
    \calE_N(t) \le \ti C \jt^{\delta} \calE_N(0),
\end{equation}
\begin{equation}\label{bdd}
    \calE_{\tilde N}(t) \le \ti C\calE_N(0),
\end{equation}
\begin{equation}\label{GE2}
    |\p_{\le N_1}| \le  \frac{\nu_N\calE_N(0)}{\la t\ra^{1/2}}, \quad |\pa\p_{\le N_1}| \le \frac{\nu_N\calE_N(0)}{\jr\mu}.
\end{equation}
for all $t \ge 0$. Recall that $\mu$ was defined in \cref*{def:mu}.
\end{theorem} 

\begin{proof} The proof will proceed via a bootstrap argument. Clearly \cref*{GE1} and \cref*{bdd} hold for small times. Assuming now that \cref*{GE1} and \cref*{bdd} hold for $0\leq t\leq T$, we improve the constant $\ti C$ by a factor of $1/2$. Thus by continuity the solution exists for all time. 

The proof that \cref*{GE1} and \cref*{bdd} hold with a better constant uses, crucially, the smallness of $|\p_{\le N_1}|$ and $|\pa\p_{\le N_1}|$. For instance, as the reader can verify below, in the region $r \le t/2$ we manage to absorb the nonlinearity to the left-hand side by way of this smallness.  We will thus prove that, for a suitable $0<\eps \ll 1$, we have the pointwise estimates 
\begin{equation}\label{GE2eps}
    |\p_{\le N_1}| \le  \frac{\eps}{\la t\ra^{1/2}}, \quad |\pa\p_{\le N_1}| \le \frac{\eps}{\jr\mu}
\end{equation}
on $0\leq t\leq T$.

Under the assumption \cref*{GE2eps}, we will show that there are constants $C_N$, $\ti C_N$ and $C_{\ti N}$, independent of $\epsilon$ and $T$, so that, for all $0\leq t\leq T$, we have 
\begin{equation}\label{E_N bound}
    \calE_N(t) \leq C_N\calE_N(0) \jt^{\ti C_N\eps},
\end{equation}
\begin{equation}\label{Nbdd}
    \calE_{\tilde N}(t) \le C_{\ti N}\calE_N(0).
\end{equation}

This means that, by picking $\ti C = \max\{2C_N, 2C_{\ti N}\}$ and $\eps < \ti C_N^{-1}\delta$, we improved $\ti C$ by a factor of $1/2$. Moreover, using \cref*{inptdcExt} and \cref*{DerProp}, we will be able to show that
\begin{equation}\label{GE2LE}
|\p_{\le N_1}| \le  \frac{ C_{N_1}\calE_{\ti N}(T)}{\la t\ra^{1/2}}, 
\qquad |\pa\p_{\le N_1}| \le \frac{C_{N_1}\calE_{\ti N}(T)}{\jr\mu}    
\end{equation}

We now pick $\eps_0 \ll 1$ small enough so that $C_{N_1} C_{\ti N} \eps_0 = \eps$. This recovers the estimate \cref*{GE2eps}.

We now present the details. Note that by continuity there is $0\leq t\leq T$ so that \cref*{GE1} and \cref*{bdd} hold for $0\leq t\leq T$. We now show that \cref*{GE2LE} holds for $0\leq t\leq T$.

In the region $8 \leq t \leq T$, Proposition~\ref{inptdcExt} directly yields%

\[
|\p_{\le N_1}| \le \frac{\bar C_{N_1}\calE_{\ti N}(T)}{\la t\ra^{1/2}}
\]

For $0\leq t\leq 8$, \cref*{GE2LE} follows immediately from the $\R^3$ Sobolev embedding
\[
\|f\|_{L^{\infty}(\R^3)} \lesssim \|f\|_{\dot H^2(\R^3)} + \|f\|_{\dot H^1(\R^3)}
\]
applied to $f=\phi_{\leq N_1}(t, \cdot)$.

We now prove the estimate for the derivative. 

In the region $8 \leq t \leq T$, we note that, by \cref*{DerProp} and the fact that $\frac{N_1+3}2 \leq N_1$ and $N_1 + 5 \le \tilde N$:
\[
\|\pa \p_{\le N_1} \|_{L^\iy(C)}  \le \bar C_{N_1} \left(\f1\jr \f1{\mu} + \|\pa\p_{\leq N_1} \|_{L^{\infty}(C)}\right)\calE_{\ti N}(T) 
\]
Using \cref*{bdd} we see that we can absorb the second term on the right to the left, since $\bar C_{N_1} \eps_0 \leq 1/2$. We thus obtain
\[
|\pa \p_{\leq N_1}| \leq 2\bar C_{N_1} \f1{\jr\mu}  \calE_{\ti N}(T) 
\]

The desired estimate \cref*{GE2LE} now follows with $C_{N_1} = 2\bar C_{N_1}$.

For $0\leq t\leq 8$, we have that $\mu \approx 1$, and \cref*{GE2LE} follows immediately from the $\R^3$ Sobolev embedding

$$\|f\|_{L^{\iy}(R < |x| < R+1)} \ls R\inv \|f_{\le 3}\|_{\lt(R - 1 < |x| < R+2)}$$
applied to $f = \pa\phi_{\leq N_1}(t, \cdot)$.

We can now assume that \cref*{GE2eps} holds on $0 \le t\le T$. We need to show that \cref*{E_N bound} and \cref*{Nbdd} hold.

We start by proving \cref*{E_N bound}. Let 
$$\calN:= \lolt + LE^*$$
be the space in which we place the nonlinearity. We will use $\lolt$ when $r\geq t/2$, and $LE^*$ when $r\leq t/2$. We remark that for this part of the proof the null structure is irrelevant.

Given the assumption of local energy decay \cref*{eq:LED}, we have 
\begin{equation}\label{from.ass}
\calE_N(t)\leq C_N \Big( \calE_N(0) + \|Q_{\le N}\|_{\calN[0,t]} \Big).
\end{equation}

It is easy to see that \cref*{Qass} implies
\[
Q_{\le N} \in S^Z(1) (\pa\phi_{\leq N_1} ) (\pa \phi_{\leq N})
\]

We define $$S_1 := \{ (s,x) : 0 \leq s \leq t, \, |x| \leq s/2 \}$$ $$S_2 = \{ (s,x) : 0 \leq s \leq t, \, s/2\leq |x| \leq 3s/2 \}$$ $$S_3 := \{ (s,x) : 0 \leq s \leq t, \, |x| \geq 3s/2 \}$$
    
  In $S_1$ we have by \cref*{GE2} 
\begin{align}\label{hgin}
\begin{split}
\|Q_{\le N}\|_{LE^*(S_1)}&\ls \| (\pa\p_{\leq N_1}) (\pa\p_{\le N})\|_{LE^*(S_1)} \\
&\ls \eps \|\jr^{-3/2}\pa\p_{\le N}\|_{LE^*(S_1)} \\
&\ls \eps \|\jr\inv \pa\p_{\le N}\|_{L^2_{t,x}(S_1)} \\
&\ls \eps \|\p_{\le N}\|_{LE^1[0,t]}.
\end{split}
\end{align}

On the other hand, in $S_2\cup S_3$ \cref*{GE2} implies that $|\pa \phi_{\le N_1}| \lesssim \epsilon\la t\ra^{-1}$. We thus obtain
\begin{equation}\label{hgout}
\|Q_{\le N}\|_{\lolt(S_2\cup S_3)} \ls \eps \|\la s\ra^{-1} \pa\p_{\le N}\|_{\lolt(S_2\cup S_3)} \ls \epsilon\int_0^t \js\inv \calE_N(s)ds
\end{equation}

We thus obtain, by \cref*{from.ass}, \cref*{hgin}, and \cref*{hgout}:
\[
\calE_N(t) \leq C_N \Big( \calE_N(0) + C\eps\int_0^t \js\inv \calE_N(s)ds\Big)
\]

By Gronwall's inequality, 
\[
    \calE_N(t) \le C_N\calE_N(0) \exp \left(C\eps \int_0^t \js\inv ds\right) \leq C_N\calE_N(0) \jt^{\ti C_N\eps}.
\]
which finishes the proof of \cref*{E_N bound}.

We will now show that \cref*{Nbdd} holds. The bounds away from the light cone are similar, and do not rely on the null structure. We first see, similarly to \cref*{hgin}, that
\begin{equation}\label{lwin}
\|Q_{\le \ti N}\|_{LE^*(S_1)} \lesssim \eps \|\p_{\le \ti N}\|_{LE^1[0,t]}
\end{equation}

Since \cref*{GE2} implies that $|\pa \phi_{\le N_1}| \lesssim \epsilon\la t\ra^{-3/2}$ in $S_3$, we obtain
\begin{equation}\label{lwout}
\|Q_{\le \ti N}\|_{L^1L^2(S_3)} \lesssim \epsilon\int_0^t \la s\ra^{-3/2} \calE_{\ti N}(s)ds  
\ls \eps \int_0^t \js^{-3/2} \calE_N(s) ds 
\ls \eps \calE_N(0)
\end{equation}where the last bound holds by \cref*{E_N bound}.

For the bound in $S_2$, however, we proceed differently. One novel aspect of the proof is that, unlike the classical null condition, the generalized null condition is not preserved after applying the vector fields in $Z$. The key observation that allows us to close the argument is that the part that does not have any tangential derivatives decays better by a factor of $1/r$. More precisely, if $Q$ satisfies the generalized null condition
\[
Q[\pa\phi, \pa\phi] \in S^Z(1) \pa\phi \overline{\pa}\phi,
\]
we will show that
\begin{equation}
Q_{\le n} \in S^Z(1) \pa\p_{\leq\x}\overline\pa\p_{\leq\xb} + S^Z\left(\jr^{-1} \right) \pa \p_{\leq\x}\pa\p_{\leq\xb}
\label{eq:vfappliedtoQ}
\end{equation}where $\x + \xb \le n$. 

Indeed, this is clear for $r\geq 1$, and it follows by induction for $r\geq 1$ provided that we show that
\begin{equation}\label{hghnull}
[Z, \tpa] \in S^Z(1) \tpa + S^Z(r^{-1}) \pa
\end{equation}

If $Z$ is a derivative, we compute
\[
[\pa_{t}, \tpa] = 0, \quad [\pa_i, \pa_v] = [\ang_i, \pa_{ r}] \in S^Z(r^{-1}) \ang, \quad [\pa_i, \ang] \in S^Z(r^{-1}) \pa.
\]
Similarly for $\Omega$ and $S$ we have
\[
[\Omega, \pa_v] = 0, \quad [\Omega, \ang] \in S^Z(1) \ang, \quad [S, \pa_v] = \pa_v, \quad [S, \ang] \in S^Z(1) \ang.
\]
This proves \cref*{hghnull}.

 Due to \cref*{eq:vfappliedtoQ}, we have
\begin{align*}
\|Q_{\le \ti N}\|_{\lolt(S_2)} &\ls \|\pa\p_{\leq N_1}\bar\pa\p_{\leq \ti N}\|_{\lolt(S_2)} + \|\bar \pa\p_{\leq N_1}\pa\p_{\leq \ti N}\|_{\lolt(S_2)} + \|\la r\ra^{-1} \pa\phi_{\leq N_1} \pa\phi_{\leq \ti N}\|_{\lolt(S_2)} .
\end{align*}

The third term is easy to handle, since $r\sim t$ in $S_2$. Using \cref*{GE2} and Gronwall's inequality we obtain
\begin{equation}\label{third}
\|\la r\ra^{-1} \pa\phi_{\leq N_1} \pa\phi_{\leq \ti N}\|_{\lolt(S_2)} \lesssim \epsilon\int_0^t \la s\ra^{-2} \calE_{\ti N}(s)ds \ls \eps \calE_N(0).
\end{equation}

We now estimate the second term. Combined with \cref*{GE2}, \cref*{tangdef} yields
\begin{equation}\label{tgpt}
\overline{\partial} \phi_{\le N_1}(t,x)|_{S_2} \ls \f{\eps}{\la t\ra^{3/2}}.
\end{equation}

By \cref*{tgpt} we have
\begin{equation}\label{lcone1}
\|\bar \pa\p_{\leq N_1}\pa\p_{\leq \ti N}\|_{\lolt(S_2)} \ls \eps \|\la s\ra^{-3/2}\pa\p_{\le \ti N}\|_{\lolt(S_2)} \ls \eps \int_0^t \js^{-3/2} \calE_N(s) ds \lesssim \eps \calE_N(0)
\end{equation} where the last bound holds by \cref*{E_N bound}.

For the first term, we estimate using \cref*{tangdef}:
\[
 \|\pa\p_{\leq N_1}\bar\pa\p_{\leq\ti N}\|_{\lolt(S_2)}  \lesssim \|\pa\p_{\leq N_1}\la s -r \ra\jr\inv\pa\p_{\leq\ti N}\|_{\lolt(S_2)} + \|\pa\p_{\leq N_1} \jr\inv \p_{\leq\ti N+1}\|_{\lolt(S_2)}
\]
By \cref*{GE2} and \cref*{GE1}, we pointwise bound $\partial \phi_{\le N_1}$, and we obtain:
\[
\|\pa\p_{\leq N_1}\la s -r \ra\jr\inv\pa\p_{\leq\ti N}\|_{\lolt(S_2)} \lesssim \eps\int_0^t \js^{-3/2}\calE_{\ti N}(s)ds \ls \eps\calE_{\ti N}(t)
\]
where we have used that inside the set $S_2$, one has $\langle s- r \rangle \ls \langle r \rangle$. Also, in $S_2$, one has $\langle s \rangle \sim \langle r \rangle$.

For the second term, we have
\[
\|\pa\p_{\leq N_1} \jr\inv \p_{\leq\ti N+1}\|_{\lolt(S_2)} \lesssim \int_0^t \eps \|\la s -r \ra^{-1/2}\jr^{-2} \p_{\leq\ti N+1}(s, \cdot)\|_{L^2(s/2 \leq |x| \leq 3s/2)} ds
\]
Clearly by Hardy's inequality
\[
\int_0^1 \eps \|\la s -r \ra^{-1/2}\jr^{-2} \p_{\leq\ti N+1}(s, \cdot)\|_{L^2(s/2 \leq |x| \leq 3s/2)} ds
\ls \int_0^1 \eps \| \pa \p_{\leq\ti N+2} (s, \cdot)\|_{L_x^2} ds \lesssim \eps \calE_N(0)
\]

On the other hand, by \cref*{GE2}, \cref*{Hardy}, \cref*{GE1}, and Hardy's inequality we obtain
\[
\begin{split}
& \int_1^t \eps \|\la s -r \ra^{-1/2}\jr^{-2} \p_{\leq\ti N+1}(s, \cdot)\|_{L^2(s/2 \leq |x| \leq 3s/2)} ds\\ 
& \ls \int_1^t \eps \js^{-2} \|\la s -r \ra^{(-1+\delta)/2} \p_{\leq\ti N+1}(s, \cdot)\|_{L^2(s/2 \leq |x| \leq 3s/2)} ds \\
& \ls \int_1^t \eps \js^{-2} \Bigl(\js^{\delta} \|\la s-r\ra^{\f{1-\delta}2}\pa_r\p_{\leq\ti N+1}\|_{L^2(s/4 \leq |x| \leq 7s/4)}+ \la s\ra^{-\f{1+\delta}2}\|\p_{\leq \ti N+1}\|_{L^2(s\approx |x|)} \Bigr) ds \\
& \ls \int_1^t \eps \js^{-3/2+\delta} \Bigl(\|\pa\p_{\leq\ti N+1}\|_{L^2(s/4 \leq |x| \leq 7s/4)} + \|\la r\ra^{-1} \p_{\leq\ti N+1}\|_{L^2(s\approx |x|)}\Bigr) ds \\
& \ls \int_1^t \eps \js^{-3/2+\delta} \Bigl(\|\pa\p_{\leq\ti N+1}\|_{L^2(s/4 \leq |x| \leq 7s/4)} + \|\la r\ra^{-1} \p_{\leq\ti N+1}\|_{L^2} \Bigr)ds \\
& \ls \int_1^t \eps \js^{-3/2+\delta} \|\pa\p_{\leq\ti N+2}\|_{L^2} ds\ls \int_1^t \eps \js^{-3/2+\delta}\calE_{\ti N+2}(s)ds \ls \eps \calE_N(0)
\end{split}
\]

We thus obtain
\begin{equation}\label{lcone}
\|Q_{\le \ti N}\|_{L^1L^2(S_2)} \lesssim \eps (\calE_{\ti N}(t) + \calE_N(0))
\end{equation}

Putting together \cref*{lwin}, \cref*{lwout} and \cref*{lcone} we get
$$\calE_{\ti N}(t) \ls \calE_N(0) + C\eps (\calE_{\ti N}(t) + \calE_N(0))$$
and \cref*{Nbdd} follows for small enough $\eps$.\end{proof}

\section{Preliminaries to the iteration}\label{sec:preliminaries}

\subsection{Overview of the proof of pointwise bounds} \label{subsec:overview}
Our starting point is the result from equation \cref*{GE2}, which provides initial decay rates for $\phi$ and its derivatives. Our ultimate goal---see the main theorem, \cref*{thm:main}---is to prove that:
\begin{equation}\label{goal}|
\phi_{\leq m}(t,x)| \lesssim \langle v \rangle^{-1} \langle u \rangle^{-1}, \quad |\pa \phi_{\leq m}(t,x)| \lesssim \langle r \rangle^{-1} \langle u \rangle^{-2}\end{equation}
where $v = t+|x|, u = t-|x|$.

For any multiindex $J$, let $w_J := S(t,0)(\phi_J(0,\cdot),\partial_t \phi_J(0,\cdot))$ denote the solution to the free wave equation $\Box W=0$ with initial data $(W(0,\cdot), \partial_t W(0,\cdot))=(\phi_J(0,\cdot),\partial_t \phi_J(0,\cdot))$, and let $\phi_J = w_J + \psi_J$.

We make an important clarification regarding notation: the function $\psi_J$ does \textit{not} coincide with $Z^J \psi$. Although for other functions in this paper the subscript $J$ consistently denotes the operator $Z^J$, we will slightly abuse this notation for convenience when dealing specifically with the function $\psi$. However, we will still use the notation 
\[
(\psi_J)_{\leq m} := \{Z^I \psi_J\}_{|I|\leq m}.
\]

Due to Kirchhoff's formula, $w_J$ decays with the desired decay rate---see \cref*{free} below---so it is enough to establish the desired bounds for $\psi_J$. To achieve this goal, we will employ weighted $L^{\infty}$-$L^{\infty}$ estimates for the fundamental solution. Our strategy involves an iterative process that progressively improves the decay rates until we reach \cref*{goal}.

Note that $\psi_J$ satisfies
\begin{equation}\label{psi}
\Box\psi_J = G_J, \quad \psi_J(0,\cdot) = (\pa_t \psi_J)(0,\cdot) =0
\end{equation}
with $G_J$ described in \cref*{first write}.

We derive improved decay for $\psi_J$ using the fundamental solution for the free wave equation in $\R^{1+3}$ and the decay rates obtained previously. In the process, we need to also obtain improved decay rates for derivatives, which is accomplished in \cref*{ss:derivative_bounds}.

From this point on, we take $n$ to be a large constant that does not depend on $m$. We allow $n$ to grow between successive estimates, but will not keep track of its precise value.

The iteration procedure is described as follows. As previously noted, $w_J$ already satisfies the required bound $\langle v\rangle^{-1}\langle u\rangle^{-1}$, and thus it suffices to focus the analysis on $\psi_J$. We begin by assuming a bound of the form
\begin{equation}\label{tembd}
|\phi_{\leq m+n}|\lesssim \langle t+r \rangle^{-\beta} \langle t-r\rangle^{-\eta}, \quad |\pa\phi_{\leq m+n}|\lesssim \la r\ra^{-\tilde\alpha} \langle t+r \rangle^{-\tilde\beta} \langle t-r\rangle^{-\tilde\eta} 
\end{equation}
where $\tilde\alpha+\tilde\beta+\tilde\eta = 1+\beta+\eta$. We then apply estimates on the fundamental solution of $\Box$, more precisely Lemmas~\ref{conversion} and ~\ref{lem:minkowski_decay}, to obtain a decay rate of the form
\[
|\psi_J|\lesssim \la r\ra^{-1} \la t+r\ra^{-\beta'} \la t-r\ra^{-\eta'}
\]
where we have $1+\beta' \geq \beta$, $\eta' \geq \eta$, and $1+\beta' + \eta' > \beta+\eta$. This already gives an improvement in the region $r\geq t/2$, but not necessarily when $r$ is small. We then apply Proposition~\ref{1st.convrsn} to turn the factor $\la r\ra^{-1}$ in \cref*{tembd} to a factor of $\la t\ra^{-1}$ in the region $r\leq t/2$. Finally, we apply Proposition~\ref{derbound} to improve the decay of derivatives.  We thus obtain
\[
|\phi_{\leq m+n}|\lesssim \langle t+r \rangle^{-1-\beta'} \langle t-r\rangle^{-\eta'}, \quad |\pa\phi_{\leq m+n}|\lesssim \la r\ra^{-\tilde\alpha'} \langle t+r \rangle^{-\tilde\beta'} \langle t-r\rangle^{-\tilde\eta'} 
\]
which is an improvement over \cref*{tembd}. We continue iterating until we obtain the bounds \cref*{goal}. At this point Lemmas~\ref{conversion} and ~\ref{lem:minkowski_decay} will not improve the decay rate, and the iteration stops.

\noindent \textbf{Note on the modified value of $\sigma$.} In the following subsections, let $\sigma$ represent the rate at which the variable coefficients approach flatness (see \cref*{first write} and the first appearance of $\sigma$ in \cref*{coeff.assu}). To simplify the iteration and avoid logarithms, we will reduce $\sigma$ as needed to a positive irrational number less than the original $\xs$, assuming $0< \xs \ll 1$. This adjustment will not affect the final decay rate.

\subsection{The pointwise estimates for the initial data}

The goal of this section is to prove the optimal decay rates for $\phi$ and its derivatives when $t=0$. More precisely, we will show that, for all $m\leq N-4$, we have

\begin{equation}\label{indata} |\phi_{\leq m}(0, x)| \ls \la r\ra^{-2}
\end{equation}
\begin{equation}\label{indatader} |\pa  \phi_{\leq m}(0, x)| \ls \la r\ra^{-3}.
\end{equation}
Because each individual summand involved is nonnegative, it is equivalent to showing that
 \begin{equation}\label{p01}
\sum_{J : |J| \le m}|(Z|_{t=0})^J(\phi_0)| +  \la r\ra \left(|\pa_x(Z|_{t=0})^J(\phi_0)| + |(Z|_{t=0})^J(\phi_1)|\right) \lesssim \la r\ra^{-2} 
 \end{equation}
In the region $r\leq 2$, \cref*{p01} follows immediately from the standard Sobolev embeddings and \cref*{eq:maintheoreminitialdata}. 

 Now fix any $R\geq 1$. Using the one-dimensional radial Poincaré inequality \cref*{eq:1dest} and the Sobolev embedding inequality on the two-dimensional sphere, we obtain for any function $f$ that 
\begin{equation}\label{SE3}
\|f\|^2_{L^\infty(A_R)} 
  \lesssim  R^{-3}\Bigl(\|f_{\leq 2}\|^2_{L^2(A_R)} 
   + R^2\|\pa_r (f_{\leq 2})\|^2_{L^2(A_R)}\Bigr).
\end{equation}

Applying \cref*{SE3} to $(Z|_{t=0})^J(\phi_0)$ we obtain
\begin{align*}
&\sum_{J:|J|\le m} \bigl\lVert (Z|_{t=0})^J(\phi_0)\bigr\rVert_{L^\infty(A_R)}^2
\lesssim
\sum_{J:|J|\le m+2} R^{-3}
\Bigl(
  \bigl\lVert (Z|_{t=0})^J(\phi_0)\bigr\rVert_{L^2(A_R)}^2
 \\& + R^2 \bigl\lVert \partial_r\bigl[(Z|_{t=0})^J(\phi_0)\bigr]\bigr\rVert_{L^2(A_R)}^2
\Bigr)
\lesssim
R^{-4}
\sum_{|J|=0}^N
\bigl\lVert (1+|\cdot|)^{1/2}\,(\phi_0)_J(\cdot)\bigr\rVert_{L^2(\mathbb{R}^3)}^2
\end{align*}
and thus
\[
\sum_{J : |J| \le m}\|(Z|_{t=0})^J(\phi_0)\|_{L^\infty(A_R)} \lesssim R^{-2}.
\]

Similarly, \cref*{SE3} applied to $(Z|_{t=0})^J(\phi_1)$ yields
\begin{align*}
& \sum_{J:|J|\le m} \bigl\lVert (Z|_{t=0})^J(\phi_1)\bigr\rVert_{L^\infty(A_R)}^2
\lesssim
\sum_{J:|J|\le m+2} R^{-3}
\Bigl(
  \bigl\lVert (Z|_{t=0})^J(\phi_1)\bigr\rVert_{L^2(A_R)}^2 
  + 
  \\ & R^2 \bigl\lVert \partial_r\bigl[(Z|_{t=0})^J(\phi_1)\bigr]\bigr\rVert_{L^2(A_R)}^2
\Bigr)\lesssim
R^{-6}
\sum_{|I|=0}^N
\bigl\lVert (1+|\cdot|)^{3/2}\,(\phi_1)_J(\cdot)\bigr\rVert_{L^2(\mathbb{R}^3)}^2
\end{align*}
and thus
\[
\sum_{J : |J| \le m}\|(Z|_{t=0})^J(\phi_1)\|_{L^\infty(A_R)} \lesssim R^{-3}.
\]

Finally, since
\[
\pa_r  = r^{-1}S|_{t=0}, \quad \ang_i = -\sum_{j\neq i}\frac{x_j}{r^{2}}\Omega_{ij}
\]
we have that
\begin{equation}\label{dert0}
\pa_x f\in S^Z(r^{-1})\sum_{|J|=1} (Z|_{t=0})^J f.
\end{equation} 

Applying \cref*{SE3} to $\pa_x (Z|_{t=0})^J(\phi_0)$ yields
\begin{align*}
\sum_{J:|J|\le m} \bigl\|\partial_x (Z|_{t=0})^J(\phi_0)\bigr\|_{L^\infty(A_R)}^2
&\lesssim
\sum_{J:|J|\le m+2} R^{-3}
  \Bigl[
    \bigl\|\partial_x (Z|_{t=0})^J(\phi_0)\bigr\|_{L^2(A_R)}^2 \\[-0.2em]
&\quad\quad
    +\,R^2\,\bigl\|\partial_r\partial_x\bigl[(Z|_{t=0})^J(\phi_0)\bigr]\bigr\|_{L^2(A_R)}^2
  \Bigr]\\[0.6em]
&\lesssim
\sum_{J:|J|\le m+3} R^{-5}
  \Bigl[
    \bigl\| (Z|_{t=0})^J(\phi_0)\bigr\|_{L^2(A_R)}^2 \\[-0.2em]
&\quad\quad
    +\,R^2\,\bigl\|\partial_r\bigl[(Z|_{t=0})^J(\phi_0)\bigr]\bigr\|_{L^2(A_R)}^2
  \Bigr]\\[0.6em]
&\lesssim
R^{-6}
\sum_{|I|=0}^N
\bigl\|(1+|\cdot|)^{\tfrac12}\,(\phi_0)_I(\cdot)\bigr\|_{L^2(\mathbb{R}^3)}^2
\end{align*}
and thus
\[
\sum_{J : |J| \le m}\|\pa_x (Z|_{t=0})^J(\phi_0)\|_{L^\infty(A_R)} \lesssim R^{-3}.
\]

\subsection{The estimates for $w$}

The goal of this section is to prove the optimal decay rates for $w$ and its derivatives. More precisely, we will show that, for $|J|+n\leq N-4$, we have

\begin{equation}\label{free} (w_J)_{\leq n}(t,x) \ls \jv\inv\ju\inv.
\end{equation}
\begin{equation}\label{freeder} \pa (w_J)_{\leq n}(t,x) \ls \jv\inv\ju^{-2}.
\end{equation}

It is enough to establish the estimates for $n=0$, since $Z^I w_J$ also solves the homogeneous wave equation.

We use Kirchhoff's formula to write, for any $x \in \R^3$ and any $t > 0$:
\begin{equation}\label{Kir}
w_J(t,x) = \f1{|\pa B(x,t)|} \int_{\pa B(x,t)} \p_J(0, y) + \pa_y \p_J(0, y) \cdot (y-x) +t \pa_t \p_J(0, y) \, dS(y).
\end{equation}
We now claim that \cref*{Kir} and \cref*{p01} yield \cref*{free} and \cref*{freeder}.

Let us start with \cref*{free}. Using \cref*{p01},  we see that, on $\pa B(x,t)$, we have
\[
|\p_J(0, y)| + |\pa_y \p_J(0, y) \cdot (y-x)| +t |\pa_t \p_J(0, y)| \lesssim \max\{\la |y| \ra^{-2}, t\la |y| \ra^{-3}\}
\]
Since $|\partial B(x,t)| = 4\pi t^2$, it thus suffices to show that, if a function $f : \R^3 \to \R$ satisfies 
\begin{equation}\label{eq:f_hypothesis}
    |f(y)| \lesssim \max\{t^{-1}\la |y| \ra^{-2}, \la |y| \ra^{-3}\}
\end{equation}
then
\begin{equation}\label{fest}
\int_{\pa B(x,t)} |f(y)| \, dS(y) \lesssim \frac{t}{\la t-r\ra\la t+r\ra}
\end{equation}

 By the triangle inequality it follows that 
\begin{equation}
    |t-r| \leq |y| \leq t+r, \quad y \in \pa B(x,t)
\label{eq:limits_of_integration}
\end{equation}
By using the hypothesis \cref*{eq:f_hypothesis} and the coarea formula
\[
\int_{\partial B(x,t)} g(y) \, dS(y) = \int_{|t-r|}^{t+r} \int_{\partial B(x,t) \cap \partial B(0,\rho)} g(y) \, dH^1(y) \, d\rho,
\] where $dH^1$ is 1-dimensional Hausdorff measure,
we obtain
\[
\int_{\pa B(x,t)} |f(y)| dS(y) \lesssim \int_{|t-r|}^{t+r} \max\{t^{-1}\la |y| \ra^{-2}, \la |y| \ra^{-3}\} |\pa B(x,t)\cap \pa B(0, \rho)| d\rho \lesssim \int_{|t-r|}^{t+r} \la\rho\ra^{-2} d\rho
\] 
Here we used that $|\pa B(x,t)\cap \pa B(0, \rho)| \ls \min\{\rho, t\}$, which is true since the intersection of the two spheres is a circle contained in $\pa B(0,\rho)\cap \pa B(x,t)$, and the radius of this circle will always be less than or equal to the radius of the smaller sphere. This is why we take the minimum of $\rho$ and $t$.

The desired estimate \cref*{fest} now follows from the inequalities:
\begin{align*}
& \int_{|t-r|}^{t+r} \langle\rho\rangle^{-2}\, d\rho \le \int_{|t-r|}^{t+r} \rho^{-2}\, d\rho = |t-r|^{-1} - (t+r)^{-1} \le \frac{2t}{|t-r|(t+r)}, \quad 1\leq |t-r|, \\
& \int_{|t-r|}^{t+r} \langle\rho\rangle^{-2}\, d\rho \le \int_{|t-r|}^{t+r} 1\, d\rho \le 2t, \quad |t-r| < 1, \quad t+r \leq 2, \\
& \int_{|t-r|}^{t+r} \langle\rho\rangle^{-2}\, d\rho \le \int_{|t-r|}^1 1\, d\rho + \int_1^{t+r} \rho^{-2}\, d\rho \le 1 + 1, \quad |t-r| < 1, \quad t+r \geq 2.
\end{align*}

The proof of \cref*{freeder} is similar. We see that \cref*{dert0} implies that
\[
|\pa^2_x [(\phi_0)_J](x)| + |\pa_x[ (\phi_1)_J](x)| \lesssim \la |x|\ra^{-3}
\]

It thus suffices to show that, if a function $f : \R^3 \to \R$ satisfies 
\begin{equation*}
    |f(y)| \lesssim \max\{t^{-1}\la |y| \ra^{-3}, \la |y| \ra^{-4}\}
\end{equation*}
then
\begin{equation*}
\int_{\pa B(x,t)} |f(y)| dS(y) \lesssim \frac{t}{\la t-r\ra^2\la t+r\ra}
\end{equation*}
which follows in a similar way as above.

\subsection{Setting up the problem for $\psi$}

We write \cref*{psi} more explicitly as
\begin{equation}\label{eq:psi_definition}
\Box\psi_J = \Box \phi_J = (\Box \phi)_J + [\Box, Z^J]\phi = Q_J + ((\Box-P)\phi)_J + [\Box, Z^J]\phi 
\end{equation}
Recall that
\[
(\Box - P)\phi = -\pa_\alpha(h^{\alpha\beta}\pa_\beta\phi +B^{\alpha}\phi) - g^\xo \Delta_\xo \phi - (V-\pa_\alpha B^{\alpha})\phi 
\]
Using the assumptions \cref*{coeff.assu}, we have
\[
((\Box-P)\phi)_J \in \pa \Bigl(S^Z(\la r\ra^{-1-\sigma}) \pa\phi_{\leq |J|}+ S^Z_\text{der}(\jr^{-1-\xs}) \phi_{\leq |J|}\Bigr) + S^Z(\la r\ra^{-2-\sigma}) \phi_{\leq |J|+2} 
\]
Similarly, we obtain for some constant $c_J$: 
\begin{align*}
&[\Box, Z^J]\phi
 = c_J (\Box\phi)_{\le |J|-1} = c_J Q_{\le |J|-1}
  + c_J \bigl((\Box-P)\phi\bigr)_{\le |J|-1} \\[0.3em]
&\in S^Z(1)\,Q_{\le |J|-1}
  + \,\partial\Bigl(S^Z(\langle r\rangle^{-1-\sigma})\,\pa\phi_{\le|J|}+ S^Z_\text{der}(\jr^{-1-\xs}) \phi_{\leq |J|}\Bigr)
  + S^Z(\langle r\rangle^{-2-\sigma})\,\phi_{\le|J|+1}
\end{align*}
 We can thus write
\begin{equation}\label{first write}
\begin{split}
\Box\psi_J
&\in
\partial\Bigl(S^Z(\langle r\rangle^{-1-\sigma})\pa\phi_{\le m} + S^Z_\text{der}(\jr^{-1-\xs}) \phi_{\leq m}\Bigr) \\ & +S^Z(\langle r\rangle^{-2-\sigma})\,\phi_{\le m+2}
\;+\;S^Z(1)\,Q_{\le m},
\quad m = |J|.
\end{split}
\end{equation}
When we commute vector fields with the null form in \cref*{eq:problem}, we obtain more than one null form, but for the purposes of pointwise decay iteration we may treat all of these null forms as a single null form because at each stage of the iteration we have a uniform upper bound on these forms. For simplicity, we continue to denote this sum of null forms by $Q$ (by a slight abuse of notation).

When $r\leq t/2$ or $r\geq 3t/2$, we gain a factor of $1/\la r \ra$ for derivatives. However, in the intermediate region $t/2 \leq r \leq 3t/2$, we only gain a factor of $1/\langle u\rangle$, which introduces an additional difficulty. To address this issue, we observe that, for any function $f$, the following identity holds:
\begin{equation}\label{D decomp}
\pa f \in S^Z(\la r\ra^{-1}) f_{\leq 1} + S^Z(1) \pa_t f, \quad \text{ in } \{ r\geq t/2\}.
\end{equation}
This is obvious for $\pa_t$ and $\ang$, whereas for $\pa_r$ we write
\[
\pa_r = \frac{S}{r}- \frac{t}{r}\pa_t.
\]

Let $\chi_{\text{cone}}$ be a cutoff function localized to the region
\[
R_{\text{cone}} = \{(t,x) : t/2 \le |x| \le 3t/2,\; t \ge 1/2\}.
\]
More precisely, we choose $\chi_{\text{cone}}$ such that it is identically $1$ in the region
\[
\{(t,x) : t/4 \le |x| \le 3t/4,\; t \ge 1\},
\]
and is supported in $R_{\text{cone}}$. Thus $\chi_{\text{cone}}$ vanishes near $t=0$.

We now rewrite \cref*{first write} by using \cref*{D decomp}. Since $
\pa \chi_\text{cone} \in S^Z(\la r\ra^{-1})$, we can write

\begin{equation}\label{decomp}
\Box\psi_J =  G_1 + G_2 + G_3,
\end{equation}
where
\begin{equation}\label{final write}
\begin{split}
G_1 \in & S^Z(\la r\ra^{-2-\sigma}) \phi_{\leq m+2} + (1- \chi_\text{cone}) \times \left(S^Z(\la r\ra^{-1-\sigma}) \times \pa\phi_{\leq m+1}\right) \\ 
& G_2 := \partial_t \tilde{G}_2, \quad \tilde{G}_2 \in \chi_\text{cone} \times \Bigl(S^Z(\la r\ra^{-1-\sigma}) \pa\phi_{\leq m}+ S^Z_\text{der}(\jr^{-1-\xs}) \phi_{\leq m}\Bigr) \\
& G_3 = Q_{\le m}
\end{split}\end{equation}
For each $i\in \{1,2,3\}$, we define 
$$\Box \psi_{J,i}  = G_i, \qquad \psi_{J,i}(0,\cdot) = \partial_t\psi_{J,i}(0,\cdot) = 0.$$

\subsection{Estimates for the fundamental solution}

We have the following result, which is similar to previous classical results, see for instance \cite{John}, \cite{As}, \cite{STz}, \cite{Sz2}. Also related (but slightly different) versions of \cref*{eq:conversionbound2} can be found in \cite{MTT}, \cite{L}. We note that Lemma \ref{conversion} holds in greater generality in terms of the range of $\alpha, \beta, \eta$, but we will state it for only the less general range of exponents that we will need in the iteration to come later. 

\begin{lemma}\label{conversion}
Assume that $\psi : [0,\iy)\times\R^3\to\R$ solves $$\Box\psi (t,x)= g(t,x), \qquad \psi(0,\cdot) = 0, \quad (\pa_t \psi)(0,\cdot) = 0. $$ 
Define
\begin{equation}\label{hdef}
h(s, \rho) = \sum_{i=0}^2 \|\Omega^i g (s, \rho\omega)\|_{L^2(\S^2)}
\end{equation}where $\Omega^i$ denotes $i$ applications of a rotation vector field. 

Consider parameters $2<\alpha <3$, $\beta\geq 0$, $0\leq\eta \neq 1$, and assume that
\begin{equation}\label{inhombds}
h(s, \rho)\lesssim  \frac{1}{\langle \rho \rangle^\alpha \langle s + \rho \rangle^\beta \langle s - \rho \rangle^\eta}.
\end{equation}
Define
\[
\tilde\eta=\min\{1, \eta\}.
\]
We then have for $t\geq 1$ that
\begin{equation}\label{eq:conversionbound2}
    |r\psi(t, x)| \lesssim  \begin{cases}
\langle t + r \rangle^{3-(\alpha+\beta+\eta)}, & \alpha + \beta + \eta < 3, \\
\la \log(1 + t + r)\ra, & \alpha + \beta + \eta = 3, \\
\langle t - r \rangle^{3-(\alpha+\beta+\tilde\eta)}, & \alpha + \beta + \eta > 3.
\end{cases}
\end{equation}

Moreover, if $g$ is supported in $R_{\text{cone}}$ and satisfies \cref*{inhombds} in the case $\alpha <2$, $\beta= 0$, $0\leq \eta< 1$, we have that
\begin{equation}\label{eq:conversionbound3}
|r\psi(t, x)| \lesssim \langle t + r \rangle^{3-(\alpha+\eta)}.
\end{equation}
\end{lemma}

\begin{proof}

By using Sobolev embedding on the sphere, and the positivity of the fundamental solution of $\Box$, it is enough to obtain bounds for the radial function $\psi$ that solves
\begin{equation}\label{1dbox}
  \Box \psi(t, r) = h(t, r), \qquad \psi[0]=(0,0).
\end{equation}

Using D'Alembert's formula, we have, for some absolute constant $c>0$,
\[
0 \leq r\psi(t,r) \leq c \int_{D_{tr}} \rho h(s,\rho) ds d\rho,
\]
where $D_{tr}$ is the backwards light cone with vertex $(r, t)$, defined as
$$D_{tr} := \{ (\rho,s)  \in \R_+^2: -(t+r) \leq s-\rho\leq t-r, \ |t-r| \leq s+\rho \leq t+r\}.
$$
We use $s$ for the temporal variable in $D_{tr}$ and $\rho$ for the radial variable in the same set.

We will now prove \cref*{eq:conversionbound2}.

Let us start with the case $|x|\geq t$, which is easier. We write $D_{tr}$ in radial ($\rho$) and null ($u' := \rho-s$) coordinates, as ${D_{tr}} = \{r-t\leq u' \leq \rho \leq r+t\}$. 
We obtain, using the fact that $\alpha+\beta > 2$:
\begin{align*}
\int_{D_{tr}} \rho h(s,\rho) ds d\rho \lesssim \int_{r-t}^{r+t} \int_{u'}^{r+t} \jrho^{1-\alpha-\beta} d\rho \, \la u'\ra^{-\eta} du' &\lesssim \int_{r-t}^{r+t} \la u'\ra^{2-(\alpha+\beta+\eta)} du' \\
&\lesssim \begin{cases}
\langle t + r \rangle^{3-(\alpha+\beta+\eta)}, & \alpha + \beta + \eta < 3, \\
\la\log(1 + t + r)\ra, & \alpha + \beta + \eta = 3, \\
\langle t - r \rangle^{3-(\alpha+\beta+\eta)}, & \alpha + \beta + \eta > 3.
\end{cases}
\end{align*} 

The case $|x|<t$ is slightly more difficult. We write
\[
D_{tr} = D^{int}_{tr}\cup D^{ext}_{tr},
\]
where
\[
D^{int}_{tr} = D_{tr}\cap\{s\leq t\}, \quad D^{ext}_{tr} = D_{tr}\cap\{s\geq t\}
\]

In the exterior region we have that $\la\rho\ra \approx \la s+\rho\ra$. We observe that we can split $D^{ext}_{tr}$ into a rectangle and a triangle which intersect only on a subset of their boundaries. We perform a change of variables. On the rectangle, we use double null coordinates, and on the triangle, we use $(\rho,u')$ coordinates (similar to how we dealt with the triangular region in the case $|x| \ge t$). For the null coordinate in this part, we choose the convention such that it takes nonnegative values, thus $u' = \rho-s$. These coordinate choices are shown in the limits of integration in the first line below. We obtain, using $\alpha+\beta>2$ and $\eta<1$:
\[
\begin{split}
\int_{D^{ext}_{tr}} \rho h(s,\rho) ds d\rho & = \int_0^{t-r}\int_{t-r}^{t+r} \rho h(s,\rho) dv' du' + \int_{t-r}^{t+r} \int_{u'}^{t+r} \rho h(s,\rho) d\rho du' \\
& \lesssim \int_0^{t-r}\int_{t-r}^{t+r} \la v'\ra^{1-\alpha-\beta}\la u'\ra^{-\eta} dv' du' + \int_{t-r}^{t+r} \int_{u'}^{t+r} \jrho^{1-\alpha-\beta} d\rho \, \la u'\ra^{-\eta} du' \\ 
&\lesssim \int_0^{t-r} \la t-r\ra^{2-\alpha-\beta} \la u'\ra^{-\eta} du' + \int_{t-r}^{t+r} \la u'\ra^{2-\alpha-\beta-\eta}du' \\
& \lesssim \begin{cases}
\langle t + r \rangle^{3-(\alpha+\beta+\eta)}, & \alpha + \beta + \eta < 3, \\
\la\log(1 + t + r)\ra, & \alpha + \beta + \eta = 3, \\
\langle t - r \rangle^{3-(\alpha+\beta+\eta)}, & \alpha + \beta + \eta > 3.
\end{cases}
\end{split}\]

For the interior region, we will show the stronger estimate
\[
\int_{D^{int}_{tr}} \rho h(s,\rho) ds d\rho \lesssim \langle t - r \rangle^{3-(\alpha+\beta+\tilde\eta)}
\]

We partition the set $D^{int}_{tr}$ dyadically:
\[
D^{int}_{tr} = \bigcup_{1\leq R \leq t}  D_{tr}^R, \quad D_{tr}^R = D^{int}_{tr}\cap A_R.
\]

For dyadic regions where $1\leq R \leq \la t-r\ra$, we use $\la\rho\ra \sim R$, $s+\rho \geq s-\rho \sim \la t-r\ra$, and the fact that $|D_{tr}^R| \approx R^2$. We obtain

\[
\int_{D_{tr}^R} \rho h(s,\rho) ds d\rho \lesssim \int_{D_{tr}^R} \la\rho\ra^{1-\alpha} \la s+\rho\ra^{-\beta} \la s-\rho\ra^{-\eta} d\rho ds  \lesssim  R^{3-\alpha} \la t-r\ra^{-\beta-\eta}, 
\]
and after dyadic summation, using that $\alpha < 3$, we obtain
\begin{equation}\label{ptwsecpt}
\sum_{R \text{ dyadic}, \, R < \la t-r\ra} \int_{D_{tr}^R} \rho H ds d\rho \lesssim  \la t-r\ra^{3-\alpha} {\la t-r\ra^{-\beta-\eta}} \lesssim \la t-r\ra^{3-(\alpha+\beta+\eta)}
\end{equation}

For dyadic regions where $\la t-r\ra \leq R\leq t$, we use that $R \sim \rho \le s+\rho $, and integrate with respect to $u'=s - \rho$ and $\rho$. We obtain
\[\begin{split}
\int_{D_{tr}^R} \rho h(s,\rho) ds d\rho & \lesssim \int_{D_{tr}^R} \la\rho\ra^{1-\alpha} \la s+\rho\ra^{-\beta} \la s-\rho\ra^{-\eta} d\rho ds  \lesssim R^{1-\alpha-\beta} \int_{0}^{t-r} \la u'\ra^{-\eta} du' \\ & \lesssim R^{2-\alpha-\beta} \la t-r\ra^{1-\tilde\eta}.
\end{split}
\]

Since $\alpha+\beta>2$, we obtain after dyadic summation 
\begin{equation}\label{ptwsefar}
\sum_{R \text{ dyadic},\, \la t-r\ra \leq R\leq t} \int_{D_{tr}^R} \rho H ds d\rho \lesssim \la t-r\ra^{3-\alpha-\beta-\tilde\eta}.
\end{equation}

The conclusion for the interior region follows from \cref*{ptwsecpt} and \cref*{ptwsefar}.

The proof for \cref*{eq:conversionbound3} is much simpler. We have
\[
\int_{D_{tr}} \rho h(s,\rho) ds d\rho \lesssim \int_{-(t+r)}^{t-r} \int_{|t-r|}^{r+t} \la v\ra^{1-\alpha} dv \, \la u\ra^{-\eta} du \lesssim \la t+r\ra^{2-\alpha} \la t+r\ra^{1-\eta} = \la t+r\ra^{3-\alpha-\eta}.
\]Here, $v := \rho+s$ and $u := \rho-s$ are temporary notations for the calculation. 
This used that $\rho \sim v$ in $R_{\text{cone}}$. 
\end{proof}

In light of the decomposition of the inhomogeneity in \cref*{decomp}, we also require the following result for an inhomogeneity of the form $\partial_t g$, supported near the cone. This result is similar to \cref*{conversion}, but the solution manifests faster decay---more precisely by an additional factor of $\la t - |x| \ra$ in the estimate.

\begin{lemma}\label{lem:minkowski_decay}
Let $\psi : [0,\infty) \times \mathbb{R}^3 \to \mathbb{R}$ solve
\begin{equation*}
\Box \psi = \partial_t g, \quad \psi(0,\cdot) = 0, \quad (\partial_t \psi)(0,\cdot) = 0,
\end{equation*}
where $g: [0,\infty) \times \mathbb{R}^3 \to \mathbb{R}$ is a given function supported in $R_{\text{cone}} = \{(s,x) : s/2 \le |x| \le 3s/2, s \ge 1/2\}$. Define, for $s \ge 1/2$ and $r' > 0$,
\begin{multline*}
h_g^*(s,r') = \sum_{i=0}^2 \left( \|\Omega^i g(s,r'\cdot)\|_{L^2(\mathbb{S}^2)} + \|\Omega^i (Sg)(s,r'\cdot)\|_{L^2(\mathbb{S}^2)} \right. \\ \left. + \|\Omega^i (\Omega g)(s,r'\cdot)\|_{L^2(\mathbb{S}^2)} + \langle s-r' \rangle \|\Omega^i (\partial g)(s,r'\cdot)\|_{L^2(\mathbb{S}^2)} \right).
\end{multline*}
Consider parameters $2 < \alpha < 3$, $0\leq\eta <1$ with $\alpha+\eta > 3$, such that
\begin{equation}\label{eq:g_bounds}
h_g^* (t,r) \ls \frac{1}{\langle r\rangle^{\alpha}\langle t-r\rangle^{\eta}}.
\end{equation}
Then, for $r=|x|$, the solution $\psi$ satisfies
\begin{equation}\label{eq:conversionbound1der}
|r\psi(t, x)| \lesssim  \langle t - r \rangle^{2-(\alpha+\eta)}.
\end{equation}

Moreover, if $g$  satisfies \cref*{eq:g_bounds} with $\alpha + \eta <3$ and  $0\leq \eta< 1$, we have that
\begin{equation}\label{eq:conversionbound2der}
|r\psi(t, x)| \lesssim \la t-r\ra^{-1} \langle t + r \rangle^{3-(\alpha+\eta)}.
\end{equation}
\end{lemma}

\begin{proof}
Let $\tilde{\psi}$ be the solution to $\Box\tilde{\psi} = g$ with zero initial data: $\tilde{\psi}(0,\cdot) = \partial_t\tilde{\psi}(0,\cdot) = 0$. The support condition on $g$ (vanishing for $s < 1/2$) ensures $g(0,\cdot)=0$. Since $[\Box, \partial_t]=0$, $\psi = \partial_t\tilde{\psi}$ is the unique solution to $\Box\psi = \partial_t g$ with the specified zero initial data.

The Lorentz boosts are defined by $L^j := t\partial_j + x^j\partial_t$ for $j \in \{1,2,3\}$. We have
$$L^j - \frac{x^j}{t}S = \left( \frac{t}{r^2} \sum_{k=1}^3 x^k \Omega_{kj} \right) + \frac{x^j(t-r)(t+r)}{rt}\partial_r.$$
Rearranging the terms,
\begin{equation}\label{eq:Lorentz_boost-identity}
L^j = \frac{x^j}{t}S + \frac{t}{r^2} \sum_{k=1}^3 x^k \Omega_{kj} + \frac{x^j(t-r)(t+r)}{rt}\partial_r.\end{equation}

Let $u = t-r$ and $v = t+r$ denote the null coordinates. The identity $uv\partial_t = tS - \sum_{i=1}^3 x^i L^i$ implies that for any smooth function $f$,
\begin{equation}\label{eq:time_deriv_est}
\la u\ra|\partial_t f| \leq |Sf| + \sum_{i=1}^3 |L^i f| + |\pa_t f|.
\end{equation}The function $|\partial_t f|$ on the right-hand side is included to make \cref*{eq:time_deriv_est} true for $-1 \leq u \leq 1$. For sufficiently large $|u|$ (e.g., $\la u \ra > 2$), the final term can be absorbed into the left-hand side, yielding $(\la u \ra/C')|\partial_t f| \le |Sf| + \sum |L^j f|$.

Let $Z \in \{S, L^1, L^2, L^3\}$. The function $\phi_{Z} := Z\tilde{\psi}$ solves $\Box \phi_{Z} = Zg + [Z,\Box]\tilde{\psi}$.
For $Z=L^j$, $[ \Box, L^j]=0$, so $\Box (L^j\tilde{\psi}) = L^j g$.
For $Z=S$, $[\Box, S]=2\Box$, so $\Box (S\tilde{\psi}) = S(\Box\tilde{\psi}) + [\Box,S]\tilde{\psi} = Sg + 2\Box\tilde{\psi} = Sg + 2g$.

Given the support of $g$, the identity \cref*{eq:Lorentz_boost-identity} and the bound \cref*{eq:g_bounds}, we can deduce that
\[ \sum_{k=0}^2 \|\Omega^k (Zg)(s,r'\cdot)\|_{L^2(\mathbb{S}^2)} \lesssim  \langle r' \rangle^{-\alpha} \langle s-r' \rangle^{-\eta} \quad \text{for } Z \in \{S, L^j\}. \]
The term $2g$ in the source for $S\tilde{\psi}$ is also controlled by $C_g \langle r' \rangle^{-\alpha} \langle s-r' \rangle^{-\eta}$ due to its explicit inclusion in the definition of $h_g^*$.

Moreover, $\tilde{\psi}\equiv 0$ for $t\leq 1/2$ due to the support properties of $g$. Applying \cref*{eq:time_deriv_est} to $\tilde{\psi}$, we obtain
\begin{equation}\label{eq:psi_est}
\la u\ra|\psi| \leq |S\tilde{\psi}| + \sum_{i=1}^3 |L^i \tilde{\psi}|+|\pa_t\tilde{\psi}|.
\end{equation}

We apply \cref*{eq:conversionbound2} to $\phi_Z$ for $Z \in \{S, L^j\}$. The source parameters are $(\alpha_{src}, \beta_{src}, \eta_{src}) = (\alpha, 0, \eta)$. Given $2 < \alpha < 3$, $0 \le \eta < 1$, and the important condition $\alpha+\eta > 3$, we are in the third case of Lemma~\ref{conversion}. 
This yields:
\[ |r (Z\tilde{\psi})(t,x)| \lesssim \langle t-r \rangle^{3-(\alpha+0+\eta)} =  \langle t-r \rangle^{3-(\alpha+\eta)}. \]

Therefore,
\begin{equation} \label{eq:Zpsi_bound}
    |(S\tilde{\psi})(t,x)| + \sum_{j=1}^3 |(L^j \tilde{\psi})(t,x)| \lesssim r^{-1} \langle t-r \rangle^{3-(\alpha+\eta)}.
\end{equation}
Separately, $\psi = \partial_t\tilde{\psi}$ solves $\Box \psi = \partial_s g$. The source $\partial_s g = \frac{1}{s}S g - \frac{r'}{s}\partial_{r'} g$. On the support $R_{\text{cone}}$ ($s \sim r'$), this source $\partial_s g$ effectively has an additional $(r')^{-1}$ spatial decay compared to $Sg$. Applying the same \cref*{eq:conversionbound2}, we obtain bounds for $|\partial_t \tilde{\psi}|$ that are at least as good as \cref*{eq:Zpsi_bound}. Thus \begin{equation} \label{eq:RHS_sum_bound}
|(S\tilde{\psi})(t,x)| + \sum_{j=1}^3 |(L^j \tilde{\psi})(t,x)| + |\psi(t,x)| \lesssim  r^{-1} \langle t-r \rangle^{3-(\alpha+\eta)}.
\end{equation}
By \cref*{eq:psi_est}, we obtain
\[ \langle t-r \rangle |\psi(t,x)| \lesssim r^{-1} \langle t-r \rangle^{3-(\alpha+\eta)}, \]
which is \cref*{eq:conversionbound1der}. \cref*{eq:conversionbound2der} follows in a similar way, using  \cref*{eq:conversionbound2} in the case $\alpha+\eta < 3$.

\end{proof}

\subsection{Improved decay in the exterior region}

We note that \cref*{GE2} gives in the region $r > \frac{3t}2$ (recall the definition of $\mu$ in \cref*{def:mu}):
\begin{equation}\label{initialext_orig} 
|\p_{\le m+n}| \ls \frac{1}{\la t\ra^{1/2}}, \quad |\pa\p_{\le m+n}| \ls \frac{1}{\la r\ra \la t\ra^{1/2}}.
\end{equation}
Here, $m+n$ denotes an appropriate high order of $Z$-derivatives, up to $N_1 = N/2$, consistent with the bounds from the global existence theorem.

This is not good enough to start the iteration process, due to the presence of factors of $t$ in the denominator. Moreover, the bound for $|\pa\p_{\le m+n}|$ is not quite strong enough to use in the following subsection. We will thus improve \cref*{initialext_orig} by replacing the $t$ factors in the denominator with $v$ factors.

First, observe that \cref*{R-Sob2} and \cref*{bdd} directly yield the inequality%
\footnote{Here, and in the subsequent application of \cref*{DerProp}, the use of \cref*{bdd} (controlling $\mathcal{E}_{\tilde{N}}$ where $\tilde{N}=N-2$) requires that the number of derivatives on $\phi$ in the $LE^1$ norms (e.g., $m+n+3$ or $N_1+5$) is at most $\tilde{N}$. This is satisfied if $N$ from the global existence theorem is chosen sufficiently large (e.g., $N \ge 14$ for $N_1+5 \le N-2$).}
\begin{equation}\label{initialextphi_orig} 
|\p_{\le m+n}| \lesssim \frac{1}{\langle r \rangle^{1/2}}, \quad r > \frac{3t}{2},
\end{equation}
which establishes the desired bound for $\p$.

 We will now prove that
\begin{equation}\label{initialextder_orig} 
|\pa\p_{\le N_1}| \ls \frac{1}{\jr^{3/2}}, \quad r>\frac{3t}2.
\end{equation}
Indeed, applying the same arguments from \cref*{DerProp} in the region $C^T_R$, we obtain
\[
\|\pa\p_{\le N_1} \|_{L^{\infty}(C_R^T)} \lesssim \f1{\jr^{\frac12}} \left(\f1{\jr} + \|\pa\p_{\le \lceil\frac{N_1+3}2\rceil} \|_{L^{\infty}(C_R^T)}\right)\|\p_{\le N_1+5} \|_{LE^1(C_R^T)}
\]
We can now use \cref*{bdd} (relying on $N_1+5 \le \tilde{N}$, see footnote) and absorb the second term on the right-hand side for $\mathcal{E}_N(0)$ sufficiently small, given that $N_1$ and $r$ are taken appropriately for the estimates.

Using \cref*{GE2} in the region $r \leq 3t/2$, and using \cref*{initialextphi_orig} and \cref*{initialextder_orig} in the region $r > 3t/2$ (with $m+n$ representing an index up to $N_1$), we conclude the stronger estimate
\begin{equation}\label{GE3} 
|\p_{\le m+n}| \ls \frac{1}{\la v\ra^{1/2}}, \quad |\pa\p_{\le m+n}| \ls \frac{1}{\la r\ra \la u\ra^{1/2}}.
\end{equation}

\subsection{Derivative bounds}\label{ss:derivative_bounds}

We will now derive better bounds for the derivatives; roughly speaking, the derivative gains $1/\la r\ra$ away from the cone, and $1/\la u\ra$ near the cone. 

We start with a lemma that establishes $L^2$ bounds for the derivatives. Recalling the notation in Sections 2.4 and 2.5, let $\calR \in \{C_T^R, C_T^U, C_R^T\}$ and denote by $\calD$ and $\ti\calR$ the enlargements described in \cref*{DTR}, \cref*{DTU} and \cref*{DRT}. Let $\nu:=\min(\jr, \ju)$, and note that
\begin{equation}\label{nu}
\begin{split}
\nu\sim R, \quad \text{in} \, \,D_T^R, D_R^T, \qquad \nu\sim U, \quad\text{in} \, \, D_T^U  
\end{split}
\end{equation}

We also define
\[
\ti\calR_0 = \ti\calR \cap \{t\geq 0\}.
\]

Note that $\ti\calR_0= \ti\calR$ except possibly when $\calR=C_R^T$.

\begin{lemma}\label{L2improv}
Let $\psi$ be any smooth enough function so that $\psi(0,x)=\pa_t\psi(0,x) = 0$.
We have
\begin{equation}\label{1stDeBd}
\|\pa \psi_{\leq m}\|_{\lt(\calR)} \ls \| \f{ \psi_{\le m+2}}{ \nu} \|_{\lt(\ti\calR_0)} + \|\jr (\Box\psi)_{\leq m} \|_{\lt(\ti \calR_0)}.
\end{equation}
\end{lemma}
\begin{proof}

The idea is that we can express all derivatives of a function $w$ as a combination of $\frac{Sw}{t-r}$, $\frac{1}r \Omega w$ and the Lagrangian $|\pa_t w|^2- |\pa_x w|^2$. The first two terms will have the desired bound, and we use the equation to control the last term. 

Note first that the estimate is trivial on $C_T^R$ and $C_R^T$ with $R\lesssim 1$, and on $C_T^U$ with $U\lesssim 1$, as $\nu\sim 1$ in such a region. We will thus assume that $\nu\gg 1$ on $\calR$.

We claim that there is a spherically symmetric cutoff function $\chi$ that is identically $1$ on $\calD$, supported in $\ti\calR$, and so that 
\begin{equation}\label{chiest}
|\partial^j \chi| \lesssim \nu^{-j}, \quad j\in \N
\end{equation}
\begin{equation}\label{boxchiest}
|\Box\chi| \lesssim r^{-1}\nu^{-1}
\end{equation}

Indeed, one can take $$\chi = \chi_1\Bigl(\frac{t}{T}, \frac{r}{R}\Bigr), \quad \calD=D_T^R$$ 
$$\chi = \chi_2\Bigl(\frac{t}{T}, \frac{t-r}{U}\Bigr), \quad \calD=D_T^U,$$
$$\chi = \chi_3\Bigl(\frac{t-T}{R}, \frac{r}{R}\Bigr), \quad \calD=D_R^T$$
for cutoffs adapted to suitable regions of size $1$. 

\cref*{chiest} clearly holds due to \cref*{nu}. \cref*{boxchiest} is trivial on $D_T^R$ and $D_R^T$, since $r\sim\nu$ on those regions. On the other hand, on $D_T^U$ we see that 
\[
(\partial_t+\partial_r)\chi = \frac1T (\partial_t \chi_2)
\]
and thus
\[
|\Box\chi| \lesssim |(\partial_t-\partial_r)(\partial_t+\partial_r)\chi| + \frac1r |\pa\chi| \lesssim r^{-1}\nu^{-1}.
\]
 
 Going back to the proof, it is enough to show that

\begin{equation}\label{1stDeBdstr}
\int_{\ti \calR_0} \chi |\pa \psi_{\leq m}|^2 \ls \| \f{ \psi_{\le m+2}}{ \nu} \|^2_{\lt(\ti\calR_0)} + \|\jr (\Box\psi)_{\leq m} \|^2_{\lt(\ti \calR_0)}.
\end{equation}

We start with the case $m = 0$, where we need to show that
\begin{equation}\label{1stDeBd0}
\int_{\ti \calR_0} \chi |\pa \psi|^2 \ls \| \f{ \psi_{\le 2}}{ \nu} \|^2_{\lt(\ti\calR_0)} + \|\jr (\Box\psi) \|^2_{\lt(\ti \calR_0)}.
\end{equation}

We begin by showing that
\begin{equation}\label{derLag}
|\pa \psi|^2 = |\partial_t \psi|^2 + |\partial_r \psi|^2 + |\ang \psi|^2 \lesssim \frac{|\Omega \psi|^2 + |S\psi|^2}{\nu^2} + \frac{t+r}{|t-r|}\left|(\pa_t \psi)^2 - |\pa_x \psi|^2\right| 
\end{equation}

Note first that
\begin{equation}
    \ang_i\psi \in S^Z\Bigl(\frac{1}{r}\Bigr) \Omega\psi.
\label{eq:angular}
\end{equation}
Since
\[
(S\psi)^2 + t(t-r)((\partial_r \psi)^2 - (\partial_t \psi)^2)- (t-r)^2 (\partial_r \psi)^2 = rt (\partial_t \psi + \partial_r \psi)^2 \geq 0
\]
we obtain that
\[
|\pa_r \psi|^2 \leq \frac{|S\psi|^2}{\nu^2} + \frac{t}{|t-r|}\left|(\pa_t \psi)^2 - |\pa_r \psi|^2\right| \leq \frac{|S\psi|^2}{\nu^2} + \frac{t}{|t-r|}\left|(\pa_t \psi)^2 - |\pa_x \psi|^2\right| + \frac{t}{|t-r|}|\ang \psi|^2
\]

Similarly, switching the roles of $r$ and $t$, we obtain
\[
|\pa_t \psi|^2 \leq \frac{|S\psi|^2}{\nu^2} + \frac{r}{|t-r|}\left|(\pa_t \psi)^2 - |\pa_x \psi|^2\right|
+\frac{r}{|t-r|}|\ang \psi|^2\]

Recalling that \cref*{eq:angular} holds, and using the fact that 
\[
\frac{t+r}{|t-r|}|\ang \psi|^2 \lesssim \frac{t+r}{r^2|t-r|}|\Omega \psi|^2 \lesssim \frac{|\Omega \psi|^2}{\nu^2}
\]
we obtain \cref*{derLag}.

We will now show that
\begin{equation}\label{LagP}
 \int_{\ti \calR_0} \chi \frac{t+r}{|t-r|}\left|(\pa_t \psi)^2 - |\pa_x \psi|^2\right| \lesssim \| \f{ \psi_{\le 2}}{ \nu} \|^2_{\lt(\ti\calR_0)} + \|\jr (\Box\psi) \|^2_{\lt(\ti \calR_0)}.
\end{equation}

  We have after integrating by parts
\begin{equation}\label{LagP2}
\int_{\ti \calR_0} \chi \psi (\Box\psi) = \int_{\ti \calR_0} \chi (-(\pa_t \psi)^2 + |\pa_x \psi|^2) + (\frac12 \Box \chi) \psi^2.  
\end{equation}

Note that there is no boundary term at $t=0$ since $\psi(0,x)= 0$.

When $\calR\in\{C_T^R, C_R^T\}$ we have $\nu\sim r$, so by Cauchy-Schwarz
\[
|\psi| |\Box\psi| \lesssim \nu^{-2}|\psi|^2 + \la r\ra^2 |\Box\psi|^2
\]
 The conclusion \cref*{LagP} now follows, using also \cref*{boxchiest}, since $\frac{t+r}{t-r}\sim 1$.

When $\calR= C_T^U$ on the other hand, we have the worse estimate $\frac{t+r}{t-r}\sim \frac{R}{U}$. To obtain \cref*{LagP}, we multiply \cref*{LagP2} by $\frac{R}{U}$ and use the Cauchy-Schwarz inequality combined with \cref*{boxchiest} to obtain
\[
\frac{R}{U}|\psi| |\Box\psi| \lesssim U^{-2}|\psi|^2 + R^2 |\Box\psi|^2
\]
\[
\frac{R}{U} |\Box\chi|\psi^2 \lesssim \nu^{-2} \psi^2.
\]
This finishes the proof of \cref*{LagP} in the region $C_T^U$.

The proof for the general case $m$ follows by induction on the number of time derivatives in $J$. Note first that
\[
\Box (Z^J\psi) \in S^Z(1)(\Box\psi)_{\leq |J|}.
\]
As long as $Z^J$ does not include any time derivatives $\pa_t$, we have
\[
Z^J\psi(0,x)=\pa_t Z^J\psi(0,x)=0
\]
and the conclusion \cref*{1stDeBdstr} follows immediately from \cref*{1stDeBd0} by using $Z^J \psi$ as a multiplier. 

If $Z^J$ contains exactly one time derivative $\pa_t$, \cref*{1stDeBdstr} also follows in the same way, since the boundary term at $t=0$ is still $0$.

Finally, if $Z^J = \pa_t^2 Z^{I}$, we use the equation to write
\[
Z^J\psi = \Delta_x Z^I\psi - \Box Z^I\psi
\]
and apply the induction hypothesis.
\end{proof}

The main result of this subsection is the following:
\begin{proposition}\label{derbound}
Let $\p$ solve %
\cref*{eq:problem}. Recall that $ \nu:=\min(\jr, \ju)$. Fix $m\in\N$, and assume that, for all $|J|\leq m$, we have
\begin{equation}\label{rbdsdrv} 
|\phi_{\le m+n}|+ |(w_J)_{\leq n}| + \nu |\pa (w_J)_{\leq n}|\ls \jr^{-\alpha}\la t+r\ra^{-\beta}\ju^{-\eta}.
\end{equation}
for some (large enough) $n=n(\sigma)$.
 We then have 
\begin{equation}
|\pa\phi_{\leq m}| \ls \jr^{-\alpha}\la t+r\ra^{-\beta}\ju^{-\eta} \nu^{-1}
\end{equation}
when $t\geq 8$.
\end{proposition}

\begin{proof}

Note first that, since $\phi_J=\psi_J +w_J$, it is enough to prove that 
\begin{equation}\label{derimp}
|\pa\psi_J| \ls \jr^{-\alpha}\la t+r\ra^{-\beta}\ju^{-\eta} \nu^{-1}
\end{equation}
for all $|J|\leq m$.

It suffices to show that, if $\calR \in \{C_T^R, C_R^T, C_U^T\}$, we have
\begin{equation}\label{derfarcone}
\|\pa\psi_J\|_{L^\infty(\calR)} \lesssim \|\nu^{-1}(\psi_J)_{\leq 5}\|_{L^\iy(\ti\calR_0)} +\|\nu^{-\sigma}(\phi_{\leq |J|+6})\|_{L^\iy(\ti\calR_0)} 
\end{equation}
for all $|J|\leq m+n-6$.

Indeed, \cref*{rbdsdrv} implies
\[
|(\psi_J)_{\leq n}| \ls \jr^{-\alpha}\la t+r\ra^{-\beta}\ju^{-\eta}
\]
and plugging into  \cref*{derfarcone} yields:
\[
\|\pa\psi_J \|_{L^\iy(\calR)} \ls \jr^{-\alpha}\jt^{-\beta}\ju^{-\eta} \nu^{-\sigma}
\]
for all $|J|\leq m+n-6$. We can now iterate the process $\lfloor \frac{1}{\sigma}\rfloor$ times. 

We now prove \cref*{derfarcone}. Note that, by \cref*{GE2} and \cref*{GE3}, we have that
\[
Q_{\leq m+n} \ls \left| \left(\pa\p_{\leq \f{m+n}2}\right)\left(\pa\p_{\le m+n} \right) \right| \lesssim \frac{1}{\la r\ra\nu^{1/2}} |\pa\p_{\le m+n}|. 
\]

We thus obtain for any $n$, taking \cref*{first write} into account:
\begin{equation}\label{Pbd}
(\Box\psi_J)_{\leq n} \lesssim \frac{1}{\la r\ra\nu^{1/2}} |\pa\p_{\le |J|+n}| + \jr^{-1-\xs}|\p_{\leq |J|+n+2}| 
\end{equation}

We now use \cref*{Pbd} combined with \cref*{L2improv} for $\psi_J$:
\begin{equation}\label{UsingPbd}
\begin{split}
 \|\pa (\psi_J)_{\le 3}\|_{\lt(\calR)} & \ls \| \nu^{-1}(\psi_J)_{\le 5}  \|_{\lt(\ti\calR_0)} + \|\la r\ra (\Box\psi_J)_{\leq 3}\|_{\lt(\ti\calR_0)} \\ 
&\lesssim \| \nu^{-1}(\psi_J)_{\le 5}  \|_{\lt(\ti\calR_0)}  + \| \nu^{-\sigma} \phi_{\leq |J|+5}\|_{\lt(\ti\calR_0)}  \\ &
\ls |\cal R|^{\frac12} \Bigl(\| \nu^{-1}(\psi_J)_{\le 5}  \|_{L^{\infty}(\ti\calR_0)} + \| \nu^{-\sigma} \phi_{\leq |J|+5}\|_{L^{\infty}(\ti\calR_0)}\Bigr).
\end{split}
\end{equation}

We now return to \cref*{SobEmbExt}, using \cref*{2nd'sEst} and \cref*{Pbd} to bound the second-order derivatives pointwise and \cref*{UsingPbd} to bound the first-order derivatives in $\lt$.
 We find
\begin{align*}
\|\pa\psi_J\|_{L^\iy(\calR)} 
&\ls |\calR|^{-\f12} \|(\pa \psi_J)_{\le 3}\|_{\lt(\calR)} + |\calR|^{-\f12} \|\jr(\Box\psi_J)_{\le 4} \|_{\lt(\calR)} \\
&\ls \| \nu^{-1}(\psi_J)_{\le 5}  \|_{L^{\infty}(\ti\calR_0)} + \| \nu^{-\sigma} \phi_{\leq |J|+6}\|_{L^{\infty}(\ti\calR_0)}.
\end{align*}

This finishes the proof of \cref*{derfarcone}.
\end{proof}

\subsection{Converting $r$ decay to $t$ decay}

In the interior region the iteration is similarly based on \cref*{conversion}, with an additional twist. It turns out that plugging in a bound of $\jr^{-1}\ju^{-\eta}$ on the right hand side is not enough to gain decay. Instead, we first need to turn the $r^{-1}$ factor into a $t^{-1}$ factor. This conversion of spatial decay to temporal decay is possible because the solution disperses and does not stay spatially concentrated for long periods of time; this leads to the fact that certain local energy decay bounds hold. These local energy decay bounds are what we will utilize here to achieve the ``conversion'' of the $r^{-1}$ factor toward a $t^{-1}$ factor.

We start with the following lemma, which is proved in Section 8 of \cite{L} (see also Proposition 3.14 in \cite{MTT}). We will reprove it below for completeness.

\begin{lemma}\label{lem:aux}

Assume that $\p$ solves \cref*{eq:problem}. We then have
\begin{equation}\label{l2int}\|\pm\|_{LE^1(\inte)} \ls T\inv \|\jr \pmn\|_{LE^1(\tinte)} + \|Q_{\le m+n}\|_{LE^*(\tinte)}.
\end{equation}
\end{lemma}
\begin{proof}
Let $\chi_\inte$ be a smooth cutoff that is identically $1$ on $\inte$, supported on $\tinte$, and so that $\pa^j \chi_\inte \lesssim T^{-j}$. Let $$\phi^{\rm{int}}:=\chi_\inte \phi.$$ Note that
\[
P\phi^{\rm{int}} = \chi_\inte Q + [P, \chi_\inte] \p
\]
By the LED estimate \cref*{eq:LED} applied to $\phi^{\rm{int}}$, we obtain
\begin{equation}\label{l2intcut}\|\phi^{\rm{int}}_{\leq m}\|_{LE^1(\inte)} \ls \|\pa\phi_{\leq m}^{\rm{int}}(\frac34T)\|_{L^2_x} + \|\left([P, \chi_\inte]\phi\right)_{\leq m}\|_{LE^*(\tinte)} + \|Q_{\le m}\|_{LE^*(\tinte)}.
\end{equation} 

Using the assumptions \cref*{coeff.assu}, and using the bounds on the derivatives of $\chi$, we easily check that
\begin{equation}\label{commcut}
\|\left([P, \chi_\inte]\phi\right)_{\leq m}\|_{LE^*(\tinte)} \lesssim T\inv \|\jr \pmn\|_{LE^1(\tinte)}.
\end{equation}
Note that the commutator $[P, \chi_\inte]\phi$ is supported in at most two dyadic annuli $A_R$, hence there is no logarithmic factor in \cref*{commcut}. 

It is thus enough to control the energy $\|\pa\phi^{\rm{int}}(T)\|_{L^2_x}$ by the space-time local energy. We will do this by averaging in time, and using the scaling vector field $S$. 

Recall that $\gamma_{T,x}$ is the integral curve of the vector field $S$ parametrized by time and such that $\gamma_{T,x}(0) = (T,x)$. By the fundamental theorem of calculus and the Cauchy--Schwarz inequality, for any $s > 0$ we have
\begin{equation}\label{Savg}
 \left| \pa u\bigl(\frac34T, x\bigr) \right|^2 \leq \frac{1}{s} \int_0^s \left| (\pa u)(\gamma_{\frac34 T,x}(\tau)) \right|^2 d\tau + \frac{s}{T^2}\int_0^s \left| (S\pa u)(\gamma_{\frac34 T,x}(\tau)) \right|^2 d\tau.
\end{equation}
uniformly with respect to $s$, so we have the freedom to choose $0 < s \leq T$ favorably depending on $x$.  

We first take $s=\frac14T$. Integrating the above estimate over $|x| \leq \frac{T}2$ yields for $u = \phi^{\rm{int}}_{\leq m}$:
\[
\int \left| \pa \phi^{\rm{int}}_{\leq m}(T) \right|^2 dx \le C T^{-1} \int_{\tinte} \left| \pa \phi^{\rm{int}}_{\leq m} \right|^2 + \left| S\pa \phi^{\rm{int}}_{\leq m} \right|^2 dx dt.
\]
Thus
\begin{equation}\label{halfgain}
\|\phi^{\rm{int}}_{\leq m}\|_{LE^1(\inte)} \ls \norm{Q_{\leq m + n}}_{LE^*(\tinte)} + T\inv \|\jr \pmn\|_{LE^1(\tinte)} + T^{-1/2} \norm{\pa \phi^{\rm{int}}_{\leq m + n}}_{L^2(\tinte)}.
\end{equation}
This takes us halfway there; indeed, we replaced $1/2$ powers of $r$ with $1/2$ powers of $t$. In order to continue, we need another Morawetz estimate:
\begin{align}\label{wwled est}
\begin{split}
\| \pa \phi^{\rm{int}}_{\leq k}\|_{L^2(\tinte)} &
\ls  \left\| \jr^{1/2} \pa \phi^{\rm{int}}_{\leq k}\Bigl(\frac34T\Bigr) \right\|_\lt  + \| \jr^{1/2} \pa \phi^{\rm{int}}_{\leq k}(T)\|_\lt \\& + \|\jr (P\phi^{\rm{int}})_{\leq k}\|_{L^2(\tinte)}, \quad k\geq 0.
\end{split}
\end{align}

We refer to Lemma 8.3 in \cite{L} for a detailed proof. For $k=0$, one uses $(r\pa_r +1)\phi^{\rm{int}}$ as a multiplier and then performs an integration by parts. One then commutes with the vector fields in $Z$. 

We now apply \cref*{wwled est} with $k=m+n$. Using \cref*{commcut}, and the fact that $\jr\leq T$ in $\tinte$, we obtain
\begin{equation}\label{Pest}
T^{-1/2} \|\jr (P\phi^{\rm{int}})_{\leq k}\|_{L^2(\tinte)} \lesssim \norm{Q_{\leq k}}_{LE^*(\tinte)} + T\inv \|\jr \p_{\leq k+n}\|_{LE^1(\tinte)}.
\end{equation}

We now apply \cref*{Savg} with $s(x)= \frac14\jr^{1/2}T^{1/2}$ and integrate in $x$. We obtain
\begin{equation}\label{12en}
\int \jr |\pa \p_{\leq k}^{\rm int} (\frac34T)|^2\,dx
    \ls \f1{{T}^{1/2}} \int_{\tinte} \jr^{1/2} |\pa \p_{\leq k}^{\rm int} |^2 + \f1{T} \jr^{3/2} |S \pa \p_{\leq k}^{\rm int} |^2 \,dxdt.
\end{equation}
The same bound holds for $\int \jr |\pa \p_{\leq k}(T)|^2\,dx$ by integrating backwards along the integral curves $\gamma$.

Since $\jr\leq T$ in $\tinte$, we obtain
\begin{equation}\label{12en2}
\begin{split}
\f1{T^{3/2}}\int_{\tinte} \jr^{3/2}|S\pa\p_{\le k}^{\rm int}|^2\,dxdt \ls \f1{T^{3/2}} \sum_{R \text{ dyadic}} R^{1/2} \int_{\tinte\cap A_R} \la r\ra |S\pa\p_{\le k}^{\rm int}|^2\,dxdt \\ \ls \f1{T^{3/2}} (\sum_{R \text{ dyadic}} R^{1/2}) \|\jr S\p_{\le k}^{\rm int}\|_{LE^1(\tinte)}^2\ls \f1{T} \|\jr S\p_{\le k}^{\rm int}\|_{LE^1(\tinte)}^2 .
\end{split}
\end{equation}
Note that there is no logarithmic loss in the summation.

Finally, by the Cauchy-Schwarz inequality we obtain for any $\delta>0$:
\begin{equation}\label{12en3}
\f1{{T}^{1/2}} \int_{\tinte} \jr^{1/2} |\pa \p_{\leq k}^{\rm int} |^2 \leq \delta \int_{\tinte} |\pa \p_{\leq k}^{\rm int} |^2 + C(\delta) \f1{T} \int_{\tinte} \jr |\pa \p_{\leq k}^{\rm int} |^2.
\end{equation}

Plugging \cref*{Pest}, \cref*{12en}, \cref*{12en2} and \cref*{12en3} into \cref*{wwled est}, we obtain after absorbing the small ($\delta\ll1$) term on the RHS of \cref*{12en3} to the LHS of \cref*{wwled est}:
\begin{equation}\label{L2final}
T^{-1/2} \| \pa \phi^{\rm{int}}_{\leq m+n}\|_{L^2(\tinte)} \lesssim \norm{Q_{\leq m + n}}_{LE^*(\tinte)} + T\inv \|\jr \pmn\|_{LE^1(\tinte)}.
\end{equation}

It is important to remark here that the index $k$ remains constant in the argument, so that the $\delta$ term can be absorbed.

The conclusion now follows from \cref*{halfgain} and \cref*{L2final}.
\end{proof}

The next proposition uses the previous lemma to turn $r$-decay in $C_T^{int}$ into $t$-decay.

\begin{proposition}\label{1st.convrsn}
Let $\phi$ solve \cref*{eq:problem}. Assume that we have in $\{r\leq 3t/4\}$:
\begin{equation}\label{rbds}
|\phi_{\le m+n}| + \jr|\pa\phi_{\le m+n}|\ls \jr^{-1}\jt^{-q}
\end{equation}
\begin{equation}\label{rbds2}
|\pa\phi_{\le m+n}| \ls \jr^{-1}\jt^{-1-q+\sigma}
\end{equation}
for some $q\geq -1/2$. We then have
\begin{equation}\label{tbds}
\|\pm\|_{L^\iy(\inte)}\ls T^{-1-q}.
\end{equation}

\end{proposition}

We explain how the next proposition will be used during the iteration. Let $k \in \mathbb{N}_{\geq 1}$ denote the current iteration stage. Hypothesis \cref*{rbds} is obtained in the first part of stage $k$, using the fundamental solution combined with \cref*{derbound}. Hypothesis
\cref*{rbds2} will hold by the bounds we obtained previously in step $k-1$.
\begin{proof}

We estimate the right hand side of \cref*{l2int}, which is $$T\inv \|\jr \pmn\|_{LE^1(\tinte)} + \|Q_{\le m+n}\|_{LE^*(\tinte)}.$$ 
Using \cref*{rbds}, we bound the first term from the RHS as follows:
\[
T^{-1} \|\jr\pmn\|_{LE^1(\tinte)} \lesssim T^{-1} T^{1/2-q} = T^{-1/2-q}.
\]

 We bound the second term of the RHS of \cref*{l2int} using \cref*{rbds2}, as follows:
\begin{align*}
\|Q_{\le m+n}\|_{LE^*(\tinte)} &\lesssim \sum_{R \text{ dyadic} : A_R \cap \inte \neq \emptyset} \|\jr^{-3/2} t^{-2-2q+2\sigma}\|_{L^2([\frac34 T, T]\times A_R)} \\
&\lesssim T^{-3/2-2q+2\sigma} \log T \lesssim T^{-1/2-q}
\end{align*}
where the last inequality holds for all $q> -1+2\sigma$. Note that unlike the estimate \cref*{commcut}, which was localized to at most two dyadic regions, this estimate involves a dyadic summation over $O(T)$ many dyadic regions, which explains the appearance of the $\log T$ factor.

Therefore, by \cref*{lem:aux}, we obtain
\[
\|\p_{\le m+n}\|_{LE^1(\inte)} \lesssim T^{-1/2-q}.
\]
The conclusion now follows, since \cref*{R-Sob} implies that
\[
\|\p_{\le m+n-3} \|_{L^\infty(C^R_T)} 
  \ls T^{-1/2}\|\p_{\le m+n}\|_{LE^1(\inte)}.
\]
\end{proof}

\section{The iteration in $\{0 \leq t \leq 8\}$}\label{sec:ext}

In this section we prove the optimal bounds \cref*{opt},\cref*{eq:derivatives} in the region $0\leq t\leq 8$. More precisely, we need to show that
\begin{equation}\label{tsmall}
|\phi_{\leq m}(t, x)|\lesssim \la r\ra^{-2}, \quad |\pa \phi_{\leq m}(t, x)|\lesssim \la r\ra^{-3}, \quad 0\leq t\leq 8.
\end{equation}

Recall that at this point we only have the bounds given by \cref*{GE2}, which are
\begin{equation}\label{tsmall1}
|\phi_{\leq m}(t, x)| \lesssim 1, \quad |\pa\phi_{\leq m}(t, x)| \lesssim \la r\ra^{-1}, \quad 0\leq t\leq 8.
\end{equation}
By using the fundamental theorem of calculus, combined with \cref*{indata} and \cref*{GE2}, we obtain
\begin{equation}\label{p01int}
|\phi_{\leq m}(t, x)| \leq |\phi_{\leq m}(0, x)| + \int_0^1 |\pa \phi_{\leq m}(s, x)| ds \lesssim \la r\ra^{-2} + \la r\ra^{-1} \lesssim \la r\ra^{-1},
\end{equation}
which improves the bound on $\phi_{\leq m}$ by a factor of $\la r\ra$. To improve the bound on the derivatives $\pa\phi_{\leq m}$, we first note that in the region $0\leq t\leq 8$ we have by Lemma~\ref{2ndDeBd'}, \cref*{p01int} and \cref*{tsmall1} that
\begin{equation}\label{2der01}
|\pa^2\phi_{\leq m}| \lesssim \la r\ra^{-1} |\phi_{\leq m+2}| + \jr^{-1}|\partial\phi_{\le m+1}| + |(\partial \phi)_{\le m}|^2 \lesssim \la r\ra^{-2},
\end{equation}
and thus the fundamental theorem of calculus, combined with \cref*{indatader} and \cref*{2der01} yields
\begin{equation}\label{dp01int}
|\pa\phi_{\leq m}(t, x)| \leq |\pa\phi_{\leq m}(0, x)| + \int_0^1 |\pa^2 \phi_{\leq m}(s, x)| ds \lesssim \la r\ra^{-3} + \la r\ra^{-2} \lesssim \la r\ra^{-2}.
\end{equation}

 The conclusion \cref*{tsmall} follows by applying the above arguments one more time, using the improved bounds 
\cref*{p01int} and \cref*{dp01int} instead of \cref*{tsmall1}.

\section{The iteration in $\{ t > 8 \}$}

In this section we prove the optimal pointwise bounds in the region $t>8$.

\begin{theorem}\label{Thmext}
Let $\phi$ be a solution to \cref*{eq:problem}. Then the following decay estimates hold in $\{ t > 8\}$:
$$|\pm| \le C_m \jv\inv \ju\inv, \quad |\pa\pm|\le C_m \jr\inv \ju^{-2}, \quad |\overline{\pa} \pm| \le C_m \jr^{-1}\jv^{-1} \ju^{-1}.$$
\end{theorem}

\begin{proof} 

We start the proof by saying that, due to Proposition~\ref{derbound}, \cref*{tangdef}, and the fact that $\nu^{-1} \approx \frac{\jv}{\jr\ju}$, we heuristically expect that
\[
\pa\pm \approx \frac{\jv}{\jr\ju} \pm, \quad \tpa{\pm}\approx \frac{\ju}{\jv}\pa\pm
\]
Our starting point is \cref*{GE3}. Taking also \cref*{tsmall} and \cref*{tangdef} into account, we have
\begin{equation}\label{1stbd}
  |\p_{\le m+n}| \ls \jv^{-1/2}, \quad |\pa\p_{\le m+n}| \ls \jr^{-1} \ju^{-1/2}, \quad \overline{\pa} \phi_{\le m+n} \ls  \jr^{-1}\jv^{-1/2}.
\end{equation}

We would now like to apply the previous results to improve the decay rates. The strategy is as follows. We apply Lemma~\ref{conversion} (and potentially Lemma~\ref{lem:minkowski_decay}) to improve the decay of $\phi_{\leq m}$ in the region $\{r\geq t/2\}$. We then apply \cref*{derbound} to improve the decay of $\pa\phi_{\leq m}$. We then apply Proposition~\ref{1st.convrsn} to improve the decay of $\phi_{\leq m}$ in the region $\{r\leq t/2\}$. Finally, we use \cref*{derbound} once again to to improve the decay of $\pa\phi_{\leq m}$ in the region $\{r\geq t/2\}$. Using the new improved bounds, we repeat the process until Lemma~\ref{conversion} stops improving the decay rate. 

This approach turns out to be infeasible if one attempts to proceed directly from the bound \cref*{1stbd}. The main obstruction arises from the inability to derive the required estimate on $\partial\phi_{\leq m}$ near the cone from the improved bounds on $\phi_{\leq m}$: Proposition~\ref{derbound} provides decay in terms of $\langle u\rangle$ near the cone, whereas decay in terms of $\langle r\rangle$ is actually required. Moreover, the presence of a factor of $\langle u\rangle^{-1}$ in the derivative’s decay rate is mildly problematic because it would introduce undesirable logarithmic terms upon application of Lemma~\ref{conversion}, which would then require additional bookkeeping.

With this in mind, we will start the iteration process with the weaker bounds
\begin{equation}\label{real1stbd}
  |\p_{\le m+n}| \ls \la v\ra^{-1/2}, \quad |\pa\p_{\le m+n}| \ls \frac{\la v\ra^{1/2}}{\la r\ra^{1+\mu} \la u\ra^{1-\mu}}.
\end{equation}
where
\[
\mu=\frac{\sigma}2
\]
Note that, except for the appearance of $\mu$, these bounds are compatible with the heuristics described at the beginning of the proof.

We will now improve the bounds \cref*{real1stbd} by $\la v\ra^{-\mu}$.

Recall the decomposition \cref*{decomp}, and define the function $H_i : [0,\infty)^2 \to [0,\infty)$ by
\[
H_i (t,r) = \sum_{k=0}^2 \|\Omega^k (G_i)_{\leq n} (t, r\cdot)\|_{L^2(\S^2)}, i\in\{1,2,3\}
\]

We obtain
\begin{equation}\label{H1est}
H_1 \lesssim \jr^{-2-\sigma}|\phi_{\leq m+n}| + (1-\chi_{cone})\jr^{-1-\sigma}|\pa\phi_{\leq m+n}|
\end{equation}
Plugging \cref*{real1stbd} into \cref*{H1est} yields:
\[
 H_1 \lesssim \la r\ra^{-2-\sigma} \la v\ra^{-1/2} + \la r\ra^{-1-\sigma}\la r\ra^{-1-\mu} \la v\ra^{-1/2+\mu} \ls \la r\ra^{-2-\sigma-\mu} \la v\ra^{-1/2+\mu}
 \]
 where we used that $u\approx v$
 on the support of $1-\chi_{cone}$.

Similarly, we also have
\begin{equation}\label{H2est}
H_2 \lesssim \chi_{cone}\jr^{-1-\sigma}|\pa\phi_{\leq m+n}| + \jr^{-2-\sigma}|\phi_{\leq m+n}|
\end{equation}
and thus
 \[
 H_2 \lesssim \la r\ra^{-1-\sigma}\frac{\la v\ra^{1/2}}{\la r\ra^{1+\mu} \la u\ra^{1-\mu}} + \la r\ra^{-2-\sigma} \la v\ra^{-1/2} \ls
 \la r\ra^{-3/2-\sigma-\mu}\ju^{-1+\mu},\]
 where we used that $r\approx v$
 on the support of $\chi_{cone}$.
 
 Finally, note that \cref*{eq:vfappliedtoQ} and  \cref*{tangdef} imply
 \begin{equation}\label{H3est}
H_3 \lesssim \frac{\ju}{\jv}|\pa\phi_{\leq m+n}|^2 + \jr^{-1}|\pa\phi_{\leq m+n}||\phi_{\leq m+n}|
 \end{equation}
 We obtain
 \[
 H_3  \ls \frac{\ju}{\jv} \frac{\la v\ra}{\la r\ra^{2+\sigma} \la u\ra^{2-\sigma}}  + \la r\ra^{-1}\frac{\la v\ra^{1/2}}{\la r\ra^{1+\mu} \la u\ra^{1-\mu}}\jv^{-1/2} \ls \frac{1}{\la r\ra^{2+\mu}\la u\ra^{1-\sigma}} \ls \frac{1}{\la r\ra^{2+\mu}\la u\ra^{1/2+\mu}},
\]
where we used $\sigma<1/3$ in the last inequality.

Note that, in the context of \cref*{conversion}, we have that $\alpha+\beta+\eta=5/2+\sigma<3$ and $\eta<1$ for all $H_i$.

By \cref*{eq:conversionbound2} with $\alpha=2+\sigma+\mu$, $\beta = 1/2-\mu$, and $\eta=0$, we obtain
\[
\psi_{J, 1} \lesssim \la v\ra^{1/2-\sigma} r^{-1}
\]
Similarly \cref*{eq:conversionbound2} with $\alpha=2+\mu$, $\beta=0$ and $\eta=1/2+\mu$ yields
\[
\psi_{J, 3} \lesssim \la v\ra^{1/2-\sigma}  r^{-1}
\]
Finally, \cref*{eq:conversionbound3} with $\alpha=3/2+\mu+\sigma$ and $\eta=1-\mu$ yields
\[
\psi_{J, 2} \lesssim \la v\ra^{1/2-\sigma}  r^{-1}.
\]

The three inequalities above, combined with the better bounds on $w$ (\cref*{free}) and \cref*{tsmall} imply the bound
\begin{equation}\label{ext1}
|\p_{\le m+n}| \ls \min\{\la v\ra^{1/2-\sigma} r^{-1}, \la v\ra^{-1/2}\} \ls \jv^{1/2-\sigma} \jr^{-1}
\end{equation}
which gains a factor of $\la v\ra^{-\sigma}$ compared to \cref*{1stbd} in the region $r\geq t/2$.

Applying \cref*{derbound} we obtain
\[
\la r\ra |\pa\p_{\le m+n}| \lesssim \jv^{1/2-\sigma} \jr^{-1}, \quad r\leq 3t/4.
\]
Recall also that \cref*{real1stbd} yields
\[
|\pa\phi_{\leq m+n}| \lesssim \jt^{\sigma/2-1/2}\jr^{-1-\sigma/2} \lesssim \la r\ra^{-1}\jt^{\sigma-1/2}, \quad r\leq 3t/4.
\]

We now use the previous bounds in \cref*{1st.convrsn} (with $q=-1/2+\sigma$) to convert $r$ decay into $t$ decay. We obtain
\begin{equation}\label{int1}
\phi_{\leq m+n} \lesssim \langle t \rangle^{-1/2-\sigma}, \quad r\leq t/2.
\end{equation}

Using \cref*{ext1} and \cref*{int1} we obtain
\[
|\p_{\le m+n}| \ls \la v\ra^{-1/2-\sigma}, 
\]
which is (better than) the desired improvement of $\la v\ra^{-\mu}$.

We now apply \cref*{derbound} and the bound above. We obtain:
\[
|\pa\p_{\le m+n}| \lesssim \frac{\jv^{1/2-\sigma}}{\jr\ju} \lesssim \frac{\jv^{1/2-\mu}}{\jr^{1+\mu}\ju^{1-\mu}}
\]
 We thus obtain the following improvement (by factors of $\la v\ra^{-\mu}$) of \cref*{real1stbd}:
\begin{equation}\label{2ndbd}
  |\p_{\le m+n}| \ls \frac{1}{\jv^{1/2+\mu}}, \quad |\pa\p_{\le m+n}| \ls \frac{\jv^{1/2-\mu}}{\jr^{1+\mu}\ju^{1-\mu}}.
\end{equation}

We now repeat the iteration, replacing \cref*{real1stbd} by \cref*{2ndbd} and applying \cref*{conversion}, followed by \cref*{1st.convrsn}, \cref*{derbound} and \cref*{tangdef}. After iterating $k=\lfloor\frac{1}{\sigma}\rfloor-1$ times, we obtain
\begin{equation}\label{3rdbd}
  |\p_{\le m+n}| \ls \frac{1}{\la v\ra^{1/2+k\mu}}, \quad \quad |\pa\p_{\le m+n}| \ls \frac{\jv^{1/2-k\mu}}{\jr^{1+\mu}\ju^{1-\mu}}.
\end{equation}

We thus have, plugging in \cref*{3rdbd} into \cref*{H1est}, \cref*{H2est} and \cref*{H3est}:
\[\begin{split}
 & H_1 \lesssim \frac{1}{\la r\ra^{2+\sigma}\la v\ra^{1/2+k\mu}}, \quad H_2 \lesssim \frac{1}{\la r\ra^{3/2+(k+1)\mu+\sigma}\la u\ra^{1-\mu}}, \\& \quad H_3  \lesssim \frac{1}{\la r\ra^{2+\sigma}\la u\ra^{1-\sigma}\jv^{2k\mu}} + \frac{1}{\la r\ra^{2+\mu}\la u\ra^{1-\mu}\jv^{2k\mu}}\ls \frac{1}{\la r\ra^{2+\sigma}\la u\ra^{1/2}\jv^{k\mu}}
\end{split}\]

Note that, in the context of \cref*{conversion}, we now have that $\alpha+\beta+\eta=5/2+(k+2)\mu > 3$ and $\eta<1$ for all $H_i$. At this point, looking at \cref*{conversion}, we switch from the case $\alpha+\beta+\eta < 3$ to the case $\alpha+\beta+\eta > 3$, and start gaining powers of $\ju$ instead of powers of $\jv$.

Indeed, applying \cref*{conversion} with $\alpha=2+\sigma$, $\beta = 1/2+k\mu$, and $\eta=0$ for $\psi_{J, 1}$;  $\alpha=3/2+(k+1)\mu+\sigma$, $\beta=0$ and $\eta=1-\mu$ for $\psi_{J, 2}$; and $\alpha=2+\sigma$, $\beta = k\mu$, and $\eta=1/2$ for $\psi_{J, 3}$, we get:
\[
|\psi_J| \lesssim r^{-1} \la u\ra^{-\tilde\mu}, \quad \tilde\mu = (k+2)\mu - 1/2 > 0,
\]

and thus as above we obtain
\begin{equation}\label{ext2}
|\p_{\le m+n}| \ls \min\{\la u\ra^{-\tilde\mu} r^{-1}, \la v\ra^{-1/2+k\mu}\} \ls \ju^{-\tilde\mu} \jr^{-1}
\end{equation}
Applying \cref*{derbound} we obtain
\[
\la r\ra |\pa\p_{\le m+n}| \lesssim \ju^{-\tilde\mu} \jr^{-1}, \quad r\leq 3t/4.
\]
Note also that \cref*{3rdbd} yields
\[
|\pa\phi_{\leq m+n}| \lesssim \la r\ra^{-1} \la u\ra^{-\tilde\mu+\sigma}
\]
Using the previous two bounds in \cref*{1st.convrsn} (with $q=\tilde\mu$) to convert $r$ decay into $t$ decay,we obtain
\begin{equation}\label{int2}
|\phi_{\leq m+n}| \lesssim \langle t \rangle^{-1} \langle u \rangle^{-\tilde\mu}, \quad r\leq t/2.
\end{equation}

Applying \cref*{derbound}, we obtain the following improved bounds 
\begin{equation}\label{4thbd}
 |\p_{\le m+n}| \ls \frac{1}{\la v\ra\ju^{\tilde\mu}}, \quad |\pa\p_{\le m+n}| \ls \frac{1}{\la r\ra \la u\ra^{1+\tilde\mu}}.
\end{equation}
 
 We now repeat the iteration, but we carefully observe that this is the last step where we can improve decay for $\psi_{J, 3}$. Indeed, we estimate 
\[ 
H_3 \lesssim \frac{1}{\la r\ra^3\la u\ra^{1+2\tilde\mu}}
\]
 and \cref*{conversion} now yields
\begin{equation}\label{3final}
\psi_{J, 3}\lesssim r^{-1} \la u\ra^{-1} 
\end{equation}
which does not improve as we gain powers of $\ju$. 

On the other hand, we can continue improving the decay rates of $\psi_{J, 1}$ and $\psi_{J, 2}$, except that we need to use \cref*{lem:minkowski_decay} for the second term in \cref*{decomp}. Define
\begin{multline*}
\tilde H_2(t,r) = \sum_{i=0}^2 \left( \|\Omega^i \tilde G_2(t,r\cdot)\|_{L^2(\mathbb{S}^2)} + \|\Omega^i (S\tilde G_2)(t,r\cdot)\|_{L^2(\mathbb{S}^2)} \right. \\ \left. + \|\Omega^i (\Omega \tilde G_2)(t,r\cdot)\|_{L^2(\mathbb{S}^2)} + \langle t-r \rangle \|\Omega^i (\partial \tilde G_2)(t,r\cdot)\|_{L^2(\mathbb{S}^2)} \right).
\end{multline*}

We estimate (see \cref*{final write}):
\begin{equation}\label{tildeH2est}
\tilde H_2(t,r) \ls \ju \la r\ra^{-1-\sigma} \pa\phi_{\leq m+n}+ \ju\jr^{-2-\xs} \phi_{\leq m+n}.
\end{equation}
As a rough guideline, we note that the first term is the larger one, because the derivative is only better by a factor of $\ju^{-1}$ near the cone.

We have, using \cref*{4thbd}:
\[
 H_1 \lesssim \frac{1}{\la r\ra^{2+\sigma}\la v\ra\ju^{\tilde\mu}}, \quad \tilde H_2 \lesssim \frac{1}{\la r\ra^{2+\sigma}\la u\ra^{\tilde\mu}}.
\]

 We now apply \cref*{conversion} for $\psi_{J, 1}$ and \cref*{lem:minkowski_decay}, more precisely \cref*{eq:conversionbound2der}, for $\psi_{J, 2}$. Note here that $2+\sigma+\tilde\mu < 3$, so unfortunately we cannot apply \cref*{eq:g_bounds}. We obtain
\[
\psi_{J, 1} \lesssim r^{-1}\ju^{-\tilde\mu-\sigma}, \quad  \psi_{J, 2} \lesssim r^{-1}\ju^{-1}\jv^{1-\sigma-\tilde\mu} 
\]

While this improves the total decay (when taking the combined sum of exponents of $\jr$, $\jv$ and $\ju$), the estimate for $\psi_{J, 2}$ is actually worse than \cref*{4thbd} near the light cone, as we had to trade decay in $v$ for decay in $u$. We thus obtain, taking \cref*{3final} and \cref*{free} into account:
\begin{equation}\label{ext3}
|\p_{\le m+n}| \ls \jr^{-1}\ju^{-1}\jv^{1-\sigma-\tilde\mu},
\end{equation}
which is the same bound as the one for $\psi_{J,2}$.

Applying \cref*{derbound} we obtain, in the set $\{r\leq 3t/4\}$:
\[
|\p_{\le m+n}| + \la r\ra |\pa\p_{\le m+n}| \lesssim \jr^{-1}\ju^{-1}\jv^{1-\sigma-\tilde\mu} \approx \jr^{-1} \jt^{-\sigma-\tilde\mu}.
\]
Recall also that the bound from the previous iteration, \cref*{4thbd}, yields in the interior region:
\[
|\pa\phi_{\leq m+n}| \lesssim \la r\ra^{-1}\jt^{-1-\tilde\mu}, \quad r\leq 3t/4.
\]

We now use the previous bounds in \cref*{1st.convrsn} (with $q=\tilde\mu+\sigma$) to convert $r$ decay into $t$ decay. We obtain
\begin{equation}\label{int3}
\phi_{\leq m+n} \lesssim \langle t \rangle^{-1-\sigma-\tilde\mu}, \quad r\leq t/2.
\end{equation}
We thus obtain the new global bound
\[
\phi_{\leq m+n} \ls \jv^{-\sigma-\tilde\mu} \ju^{-1}
\]

Applying \cref*{derbound}, we obtain the following bounds 
\begin{equation}\label{5thbd}
 |\p_{\le m+n}| \ls \jv^{-\sigma-\tilde\mu} \ju^{-1}, \quad |\pa\p_{\le m+n}| \ls \frac{\jv^{1-\sigma-\tilde\mu}}{\la r\ra \la u\ra^{2}}.
\end{equation}

 We now improve \cref*{5thbd} iteratively as before, gaining factors of $\jv^{-\mu}$. We have, using \cref*{5thbd}:
\[
 H_1 \lesssim \la r\ra^{-2-\sigma}\jv^{-\sigma-\tilde\mu} \ju^{-1+\mu}, \quad \tilde H_2 \lesssim \frac{1}{\la r\ra^{1+2\sigma+\tilde\mu}\la u\ra^{1-\mu}}.
\]
where once again we used $-1+\mu$ as a power for $\ju$ to avoid appearance of logarithms.

 We now apply \cref*{conversion} for $\psi_{J, 1}$ and \cref*{lem:minkowski_decay}, more precisely \cref*{eq:conversionbound2der}, for $\psi_{J, 2}$. We obtain
\[
\psi_{J, 1} \lesssim r^{-1}\ju^{-\tilde\mu-\sigma-\mu}, \quad  \psi_{J, 2} \lesssim r^{-1}\ju^{-1}\jv^{1-\sigma-\tilde\mu} 
\]
 which yields, taking \cref*{3final} and \cref*{free} into account:
 \[
 \phi_{\leq m+n} \ls \jr^{-1}\ju^{-1}\jv^{1-\sigma-\tilde\mu-\mu}.
 \]

 Using \cref*{derbound} and \cref*{1st.convrsn} as above we obtain
 \begin{equation}\label{6thbd}
|\p_{\le m+n}| \ls \jv^{-\sigma-\tilde\mu-\mu} \ju^{-1}, \quad |\pa\p_{\le m+n}| \ls \frac{\jv^{1-\sigma-\tilde\mu-\mu}}{\la r\ra \la u\ra^{2}}.
 \end{equation}

  We now iterate, obtaining the bounds
 \begin{equation}\label{7thbd}
|\p_{\le m+n}| \ls \jv^{-\sigma-\tilde\mu-k\mu} \ju^{-1}, \quad |\pa\p_{\le m+n}| \ls \frac{\jv^{1-\sigma-\tilde\mu-k\mu}}{\la r\ra \la u\ra^{2}}.
 \end{equation}  
Choose $k$ to be the smallest integer so that $\mu_1:=2\sigma+\tilde\mu+k\mu>1$, and let $\tilde\mu_1 := \frac{\mu_1-1}2 > 0$.

 At this point, we have
 \[
  H_1 \lesssim \la r\ra^{-2-\sigma}\jv^{-\mu_1+\sigma} \ju^{-1+\tilde\mu_1} 
 \]
 
 \[
 \tilde H_2 \lesssim \frac{1}{\la r\ra^{1+\mu_1}\la u\ra^{1-\tilde\mu_1}}.
 \]

We still apply \cref*{conversion} for $\psi_{J, 1}$. However, since $$1+\mu_1 > 2, \quad (1+\mu_1) + (1-\tilde\mu_1)=\frac{5+\mu_1}2 >3 ,$$
we can switch to \cref*{eq:conversionbound1der} for $\psi_{J, 2}$. We obtain
\[
|\psi_{J, 1}| + |\psi_{J, 2}| \lesssim r^{-1}\ju^{-1}
\]

Combining the above estimate, \cref*{3final} and \cref*{free} yields 
 \[
 \phi_{\leq m+n} \ls \jr^{-1}\ju^{-1}.
 \]

 Using \cref*{derbound} and \cref*{1st.convrsn} as above we obtain
 \begin{equation}
|\p_{\le m+n}| \ls \jv^{-1} \ju^{-1}, \quad |\pa\p_{\le m+n}| \ls \la r\ra^{-1} \la u\ra^{-2}.
 \end{equation}

 Using \cref*{tangdef} to estimate the tangential derivatives finishes the proof.
\end{proof}

\subsection*{Acknowledgements}
Part of this work was conducted while S. Looi was at UC Berkeley and he thanks UC Berkeley for its hospitality.

\end{document}